\documentclass[11pt,leqno]{article}
\usepackage{amsfonts}

\begin{document}

\def\R{{\mathbb R}}
\def\T{{\mathbb T}}
\def\S{{\mathbb S}}
\def\C{{\mathbb C}}
\def\Z{{\mathbb Z}}
\def\N{{\mathbb N}}
\def\H{{\mathbb H}}
\def\B{{\mathbb B}}
\def\diam{\mbox{\rm diam}}
\def\sn{\S^{n-1}}
\def\rr{{\cal R}}
\def\mt{{\Lambda}}
\def\e{\emptyset}
\def\dQ{\partial Q}
\def\dk{\partial K}
\def\endofproof{{\rule{6pt}{6pt}}}
\def\di{\displaystyle}
\def\dist{\mbox{\rm dist}}
\def\sa+{\Sigma_A^+}
\def\du{\frac{\partial}{\partial u}}
\def\dv{\frac{\partial}{\partial v}}
\def\dt{\frac{d}{d t}}
\def\dx{\frac{\partial}{\partial x}}
\def\con{\mbox{\rm const }}
\def\nn{{\cal N}}
\def\mm{{\cal M}}
\def\kk{{\cal K}}
\def\ll{{\cal L}}
\def\vv{{\cal V}}
\def\bb{{\cal B}}
\def\ma{\mm_{a}}
\def\lab{L_{ab}}
\def\mabn{\mm_{a}^N}
\def\man{\mm_a^N}
\def\labn{L_{ab}^N}
\def\fa{f^{(a)}}
\def\ff{{\cal F}}
\def\i{{\bf i}}
\def\gge{{\cal G}_\epsilon}
\def\gej{\chi^{(j)}_\mu}
\def\ge{\chi_\epsilon}
\def\geo{\chi^{(1)}_\mu}
\def\get{\chi^{(2)}_\mu}
\def\gei{\chi^{(i)}_{\mu}}
\def\gee{\chi_{\mu}}
\def\gett{\chi^{(2)}_{\mu}}
\def\geol{\chi^{(1)}_{\ell}}
\def\getl{\chi^{(2)}_{\ell}}
\def\geil{\chi^{(i)}_{\ell}}
\def\gee{\chi_{\ell}}
\def\tt{{\cal T}}
\def\uu{{\cal U}}
\def\wloc{W_{\epsilon}}
\def\Int{\mbox{\rm Int}}
\def\dist{\mbox{\rm dist}}
\def\pr{\mbox{\rm pr}}
\def\pp{{\cal P}}
\def\aa{{\cal A}}
\def\cc{{\cal C}}
\def\supp{\mbox{\rm supp}}
\def\Arg{\mbox{\rm Arg}}
\def\In{\mbox{\rm Int}}
\def\con{\mbox{\rm const}\;}
\def\Re{\mbox{\rm Re}}
\def\li{\mbox{\rm li}} 
\def\Seo{S^*_\epsilon(\Omega)}
\def\sdk{S^*_{\dk}(\Omega)}
\def\lae{\Lambda_{\epsilon}}
\def\ep{\epsilon}
\def\oo{{\cal O}}
\def\be{\begin{equation}}
\def\ee{\end{equation}}
\def\beqn{\begin{eqnarray*}}
\def\eeqn{\end{eqnarray*}}
\def\Pr{\mbox{\rm Pr}}

\def\gi{\gamma^{(i)}}
\def\ii{{\imath }}
\def\jj{{\jmath }}
\def\II{{\cal I}}
\def\ccij{ \cc_{i'_0,j'_0}[\eta]}
\def\dd{{\cal D}}
\def\la{\langle}
\def\ra{\rangle}
\def\bs{\bigskip}
\def\xio{\xi^{(0)}}
\def\xo{x^{(0)}}
\def\zo{z^{(0)}}
\def\Con{\mbox{\rm Const}\;}
\def\do{\partial \Omega}
\def\dk{\partial K}
\def\dl{\partial L}
\def\ll{{\cal L}}
\def\kk{{\cal K}}
\def\kk{{\cal K}}
\def\pr{{\rm pr}}
\def\ff{{\cal F}}
\def\G{{\cal G}}
\def\C{{\bf C}}
\def\dist{{\rm dist}}
\def\dds{\frac{d}{ds}}
\def\con{{\rm const}\;}
\def\Con{{\rm Const}\;}
\def\di{\displaystyle}
\def\oo{\mbox{\rm O}}
\def\hess{\mbox{\rm Hess}}
\def\gi{\gamma^{(i)}}
\def\endofproof{{\rule{6pt}{6pt}}}
\def\xm{x^{(m)}}
\def\vm{\varphi^{(m)}}
\def\km{k^{(m)}}
\def\dm{d^{(m)}}
\def\kam{\kappa^{(m)}}
\def\dem{\delta^{(m)}}
\def\xim{\xi^{(m)}}
\def\ep{\epsilon}
\def\ms{\medskip}
\def\ex{\mbox{\rm extd}}

\def\clip{C^{\mbox{\footnotesize \rm Lip}}}
\def\wlocs{W^s_{\mbox{\footnote\rm loc}}}
\def\Lip{\mbox{\rm Lip}}

\def\Xr{X^{(r)}}
\def\lip{\mbox{{\footnotesize\rm Lip}}}
\def\Vol{\mbox{\rm Vol}}

\def\naf{\nabla f(z)}
\def\so{\sigma_0}
\def\Xo{X^{(0)}}
\def\z1{z^{(1)}}
\def\Vo{V^{(0)}}
\def\Yo{Y{(0)}}

\def\uo{u^{(0)}}
\def\vo{v^{(0)}}
\def\no{\nu^{(0)}}
\def\psa{\partial^{(s)}_a}
\def\hcd{\hc^{(\delta)}}
\def\Md{M^{(\delta)}}
\def\Uo{U^{(1)}}
\def\Ut{U^{(2)}}
\def\Uj{U^{(j)}}
\def\no{n^{(1)}}
\def\nt{n^{(2)}}
\def\nj{n^{(j)}}
\def\ccm{\cc^{(m)}}

\def\ooo{\oo^{(1)}}
\def\oot{\oo^{(2)}}
\def\ooj{\oo^{(j)}}
\def\fo{f^{(1)}}
\def\ft{f^{(2)}}
\def\fj{f^{(j)}}
\def\wo{w^{(1)}}
\def\wt{w^{(2)}}
\def\wj{w^{(j)}}
\def\Vo{V^{(1)}}
\def\Vt{V^{(2)}}
\def\Vj{V^{(j)}}

\def\Ul{U^{(\ell)}}
\def\Uj{U^{(j)}}
\def\wl{w^{(\ell)}}
\def\Vl{V^{(\ell)}}
\def\Ujj{U^{(j+1)}}
\def\wjj{w^{(j+1)}}
\def\Vjj{V^{(j+1)}}
\def\Ujo{U^{(j_0)}}
\def\wjo{w^{(j_0)}}
\def\Vjo{V^{(j_0)}}
\def\vj{v^{(j)}}
\def\vl{v^{(\ell)}}

\def\f0{f^{(0)}}

\def\gl{\gamma_\ell}
\def\id{\mbox{\rm id}}
\def\piU{\pi^{(U)}}

\def\cca{C^{(a)}}
\def\bba{B^{(a)}}
\def\co{\; \stackrel{\circ}{C}}

\def\oV{\overline{V}}
\def\saa{\Sigma^+_A}
\def\sa{\Sigma_A}
\def\mta{\Lambda(A, \tau)}
\def\mtaa{\Lambda^+(A, \tau)}

\def\Int{\mbox{\rm Int}}
\def\epo{\ep^{(0)}}
\def\pH{\partial \H^{n+1}}
\def\sh{S^*(\H^{n+1})}
\def\zoo{z^{(1)}}
\def\yoo{y^{(1)}}
\def\xoo{x^{(1)}}


\def\supp{\mbox{\rm supp}}
\def\Arg{\mbox{\rm Arg}}
\def\In{\mbox{\rm Int}}
\def\diam{\mbox{\rm diam}}
\def\e{\emptyset}
\def\endofproof{{\rule{6pt}{6pt}}}
\def\di{\displaystyle}
\def\dist{\mbox{\rm dist}}
\def\con{\mbox{\rm const }}
\def\Box{\spadesuit}
\def\Int{\mbox{\rm Int}}
\def\dist{\mbox{\rm dist}}
\def\pr{\mbox{\rm pr}}
\def\be{\begin{equation}}
\def\ee{\end{equation}}
\def\beqn{\begin{eqnarray*}}
\def\eeqn{\end{eqnarray*}}
\def\la{\langle}
\def\ra{\rangle}
\def\bs{\bigskip}
\def\Con{\mbox{\rm Const}\;}
\def\clip{C^{\mbox{\footnotesize \rm Lip}}}
\def\wlocs{W^s_{\mbox{\footnote\rm loc}}}
\def\Lip{\mbox{\rm Lip}}
\def\lip{\mbox{\footnotesize\rm Lip}}
\def\Re{\mbox{\rm Re}}
\def\li{\mbox{\rm li}} 
\def\ep{\epsilon}
\def\ms{\medskip}
\def\dds{\frac{d}{ds}}
\def\oo{\mbox{\rm O}}
\def\hess{\mbox{\rm Hess}}
\def\id{\mbox{\rm id}}
\def\ii{{\imath }}
\def\jj{{\jmath }}
\def\graph{\mbox{\rm graph}}
\def\span{\mbox{\rm span}}

\def\i{{\bf i}}
\def\C{{\bf C}}

\def\ss{{\cal S}}
\def\tt{{\cal T}}
\def\E{{\cal E}}
\def\rr{{\cal R}}
\def\nn{{\cal N}}
\def\mm{{\cal M}}
\def\kk{{\cal K}}
\def\ll{{\cal L}}
\def\vv{{\cal V}}
\def\ff{{\cal F}}
\def\hh{{\cal H}}
\def\tt{{\cal T}}
\def\uu{{\cal U}}
\def\cc{{\cal C}}
\def\pp{{\cal P}}
\def\aa{{\cal A}}
\def\oo{{\cal O}}
\def\II{{\cal I}}
\def\dd{{\cal D}}
\def\ll{{\cal L}}
\def\ff{{\cal F}}
\def\G{{\cal G}}
\def\hhs{\hh^s}
\def\thhs{\widetilde{\hh}^s}
\def\hhhs{\widehat{\hh}^s}

\def\hs{\hat{s}}
\def\hz{\hat{z}}
\def\hL{\hat{L}}
\def\hl{\hat{l}}
\def\hl{\hat{l}}
\def\hc{\hat{\cc}}
\def\hbb{\widehat{\cal B}}
\def\hu{\hat{u}}
\def\hX{\hat{X}}
\def\hx{\hat{x}}
\def\hu{\hat{u}}
\def\hv{\hat{v}}
\def\hQ{\hat{Q}}
\def\hC{\widehat{C}}
\def\hF{\hat{F}}
\def\hf{\hat{f}}
\def\hii{\hat{\ii}}
\def\hr{\hat{r}}
\def\hq{\hat{q}}
\def\hy{\hat{y}}
\def\hZ{\widehat{Z}}
\def\hz{\hat{z}}
\def\hE{\widehat{E}}
\def\hR{\widehat{R}}
\def\hell{\hat{\ell}}
\def\hs{\hat{s}}
\def\hW{\widehat{W}}
\def\hS{\widehat{S}}
\def\hV{\widehat{V}}
\def\hB{\widehat{B}}
\def\hhh{\widehat{\cal H}}
\def\hK{\widehat{K}}
\def\hU{\widehat{U}}
\def\hhh{\widehat{\hh}}
\def\hdd{\widehat{\dd}}
\def\hZ{\widehat{Z}}

\def\hal{\hat{\alpha}}
\def\hbe{\hat{\beta}}
\def\hg{\hat{\gamma}}
\def\hrho{\hat{\rho}}
\def\hd{\hat{\delta}}
\def\hphi{\hat{\phi}}
\def\hmu{\hat{\mu}}
\def\hnu{\hat{\nu}}
\def\hsi{\hat{\sigma}}
\def\htau{\hat{\tau}}
\def\hpi{\hat{\pi}}
\def\hep{\hat{\epsilon}}
\def\hxi{\hat{\xi}}
\def\hLa{\widehat{\Lambda}^u}
\def\hPhi{\widehat{\Phi}}
\def\hPsi{\widehat{\Psi}}
\def\hPhii{\widehat{\Phi}^{(i)}}

\def\tc{\tilde{C}}
\def\tg{\tilde{\gamma}}  
\def\tV{\widetilde{V}}
\def\tC{\widetilde{\cc}}
\def\tr{\tilde{R}}
\def\tb{\tilde{b}}
\def\tt{\tilde{t}}
\def\tx{\tilde{x}}
\def\tp{\tilde{p}}
\def\tz{\tilde{Z}}
\def\tZ{\tilde{Z}}
\def\tF{\tilde{F}}
\def\tf{\tilde{f}}
\def\tp{\tilde{p}}
\def\te{\tilde{e}}
\def\tv{\tilde{v}}
\def\tu{\tilde{u}}
\def\tw{\tilde{w}}
\def\ts{\tilde{\sigma}}
\def\tr{\tilde{r}}
\def\tU{\tilde{U}}
\def\tS{\tilde{S}}
\def\tP{\widetilde{\Pi}}
\def\ttau{\tilde{\tau}}
\def\tLip{\widetilde{\Lip}}
\def\tz{\tilde{z}}
\def\tS{\tilde{S}}
\def\tts{\tilde{\sigma}}
\def\tVl{\widetilde{V}^{(\ell)}}
\def\tVj{\widetilde{V}^{(j)}}
\def\tVo{\widetilde{V}^{(1)}}
\def\tVj{\widetilde{V}^{(j)}}
\def\tPsi{\tilde{\Psi}}
 \def\tp{\tilde{p}}
 \def\tVjo{\widetilde{V}^{(j_0)}}
\def\tvj{\tilde{v}^{(j)}}
\def\tVjj{\widetilde{V}^{(j+1)}}
\def\tvl{\tilde{v}^{(\ell)}}
\def\tVt{\widetilde{V}^{(2)}}
\def\tR{\tilde{R}}
\def\tQ{\tilde{Q}}
\def\oL{\tilde{\Lambda}}
\def\tq{\tilde{q}}
\def\tx{\tilde{x}}
\def\ty{\tilde{y}}
\def\tz{\tilde{z}}
\def\txo{\tilde{x}^{(0)}}
\def\tso{\tilde{\sigma}_0}
\def\tmt{\tilde{\Lambda}}
\def\tg{\tilde{g}}
\def\tsi{\tilde{\sigma}}
\def\ttt{\tilde{t}}
\def\tC{\tilde{C}}
\def\tc{\tilde{c}}
\def\tell{\tilde{\ell}}
\def\trho{\tilde{\rho}}
\def\ts{\tilde{s}}
\def\tB{\widetilde{B}}
\def\thh{\widetilde{\cal H}}
\def\tV{\widetilde{V}}
\def\trr{\tilde{r}}
\def\te{\tilde{e}}
\def\tv{\tilde{v}}
\def\tu{\tilde{u}}
\def\tw{\tilde{w}}
\def\trho{\tilde{\rho}}
\def\tell{\tilde{\ell}}
\def\tz{\tilde{Z}}
\def\tF{\tilde{F}}
\def\tf{\tilde{f}}
\def\tp{\tilde{p}}
\def\ttau{\tilde{\tau}}
\def\tz{\tilde{z}}
\def\tg{\tilde{\gamma}}  
\def\tV{\widetilde{V}}
\def\tC{\widetilde{\cc}}
\def\tLa{\widetilde{\Lambda}^u}
\def\tR{\widetilde{R}}
\def\tr{\tilde{r}}
\def\tc{\widetilde{C}}
\def\tD{\widetilde{D}}
\def\tt{\tilde{t}}
\def\tx{\tilde{x}}
\def\tp{\tilde{p}}
\def\tS{\tilde{S}}
\def\tts{\tilde{\sigma}}
\def\tZ{\widetilde{Z}}
\def\tdelta{\tilde{\delta}}
\def\th{\tilde{h}}
\def\tB{\widetilde{B}}
\def\thh{\widetilde{\hh}}
\def\tep{\tilde{\ep}}
\def\tE{\widetilde{E}}
\def\tu{\tilde{u}}
\def\txi{\tilde{\xi}}
\def\teta{\tilde{\eta}}

\def\sr{{\sc r}}
\def\mt{{\Lambda}}
\def\do{\partial \Omega}
\def\dk{\partial K}
\def\dl{\partial L}
\def\wloc{W_{\epsilon}}
\def\piU{\pi^{(U)}}
\def\Rio{\R_{i_0}}
\def\Ri{\R_{i}}
\def\Rii{\R^{(i)}}
\def\Riii{\R^{(i-1)}}
\def\hRii{\widehat{\R}_i}
\def\hRiio{\widehat{\R}_{(i_0)}}
\def\Eii{E^{(i)}}
\def\Eio{E^{(i_0)}}
\def\Rj{\R_{j}}
\def\Vio{{\cal V}^{i_0}}
\def\Vi{{\cal V}^{i}}
\def\Wio{W^{i_0}}
\def\Wioo{W^{i_0-1}}
\def\hi{h^{(i)}}
\def\Psii{\Psi^{(i)}}
\def\pii{\pi^{(i)}}
\def\piii{\pi^{(i-1)}}
\def\gxyii{g_{x,y}^{i-1}}
\def\span{\mbox{\rm span}}
\def\Jac{\mbox{\rm Jac}}
\def\Vol{\mbox{\rm Vol}}
\def\limp{\lim_{p\to\infty}}
\def\hh{{\mathcal H}}

\def\yijl{Y_{i,j}^{(\ell)}}
\def\xijl{X_{i,j}^{(\ell)}}
\def\hyijl{\widehat{Y}_{i,j}^{(\ell)}}
\def\hxijl{\widehat{X}_{i,j}^{(\ell)}}
\def\eijl{\eta_{i,j}^{(\ell)}}
\def\J{\sf J}
\def\Gl{\Gamma_\ell}

\def\hLao{\widehat{\Lambda}^{u,1}}
\def\tLao{\widetilde{\Lambda}^{u,1}}
\def\Lao{\Lambda^{u,1}}
\def\cLao{\check{\Lambda}^{u,1}}
\def\cB{\check{B}}
\def\tpi{\tilde{\pi}}

\begin{center}
{\Large\bf  Pinching conditions, linearization and regularity\\ of  Axiom A flows}
\end{center}

\begin{center}
{\sc Luchezar Stoyanov}\\
{\it University of Western Australia, Perth WA 6009, Australia\\
E-mail address: stoyanov@maths.uwa.edu.au}
\end{center}

\bs


\noindent
{\it Abstract.} In this paper we study a certain regularity property of 
$C^2$ Axiom A flows $\phi_t$ over basic sets $\mt$ related to diameters of balls in Bowen's
metric, which we call regular distortion along unstable manifolds. The motivation to
investigate the latter comes from the study of spectral properties of Ruelle transfer 
operators in \cite{kn:St1}.
We prove that if the bottom of the spectrum of $d\phi_t$ over $E^u_{|\mt}$ is point-wisely pinched and
integrable, then the flow has regular distortion along unstable manifolds over $\mt$.
In the process, under the same conditions, we show that locally the flow is Lipschitz conjugate to 
its linearization over the `pinched part' of the unstable tangent bundle.

\normalsize

\section{Introduction}

\renewcommand{\theequation}{\arabic{section}.\arabic{equation}}


Let $\phi_t : M \longrightarrow M$ be a $C^2$ Axiom A flow on a
$C^2$ complete (not necessarily compact) Riemann manifold $M$ and let $\mt$ be a basic set 
for $\phi_t$.
Let $\|\cdot \|$ be the {\it norm} on $T_xM$ determined by the Riemann metric on $M$ and let
$E^u(x)$ and $E^s(x)$ ($x\in \mt$)  be the tangent spaces to the strong unstable and stable manifolds 
$W^u_\ep(x)$  and $W^s_\ep(x)$, respectively (see section 2). For any $x \in \mt$, $T > 0$ and 
$\delta\in (0,\ep]$ set
$$B^u_T (x,\delta) = \{ y\in W^u_{\ep}(x) : d(\phi_t(x), \phi_t(y)) \leq \delta \: \: , 
\:\:  0 \leq t \leq T \}\;.$$

We will say that $\phi_t$ has a {\it regular distortion along unstable manifolds} over
the basic set $\mt$  if there exists a constant $\ep_0 > 0$ with the following properties:

\ms 

(a) For any  $0 < \delta \leq   \ep \leq \ep_0$ there exists a constant $R =  R (\delta , \ep) > 0$ such that 
\be
\diam( \mt \cap B^u_T(z ,\ep))   \leq R \, \diam( \mt \cap B^u_T (z , \delta))
\ee
for any $z \in \mt$ and any $T > 0$.

\ms

(b) For any $\ep \in (0,\ep_0]$ and any $\rho \in (0,1)$ there exists $\delta  \in (0,\ep]$
such that for  any $z\in \mt$ and any $T > 0$ we have
$\diam ( \mt \cap B^u_T(z ,\delta))   \leq \rho \; \diam( \mt \cap B^u_T (z , \ep))\;.$

\bs

Part (a) of the above condition resembles the Second Volume Lemma of Bowen and Ruelle \cite{kn:BR} 
about balls in Bowen's metric; this time however we deal with diameters instead of volumes.
Given a coding of the flow $\phi_t$ over $\mt$ by means of a Markov family (see \cite{kn:B}, \cite{kn:KH} or
\cite{kn:PP}), the above properties translate into some `natural' properties of cylinders (see Sect. 3 in \cite{kn:St1}).

The aim of this paper is to describe a rather general class of flows on basic sets 
satisfying this condition.

The motivation for this work came from \cite{kn:St1}, where 
under this condition, Lipschitzness of the local stable holonomy 
maps and a certain non-integrability condition we prove strong spectral 
estimates for arbitrary potentials over basic sets for Axiom A flows, similar to 
those established by Dolgopyat \cite{kn:D}
for geodesic flows on compact surfaces (for general potentials)
and transitive Anosov flows on compact manifolds  with $C^1$ jointly non-integrable 
horocycle foliations (for the Sinai-Bowen-Ruelle potential).
It should be remarked that strong spectral estimates for Ruelle transfer operators 
lead to deep results in a variety of areas which are difficult 
(if not impossible) to obtain by  other means (see e.g. \cite{kn:PoS1}, \cite{kn:PoS2}, \cite{kn:PoS3}, \cite{kn:An}, 
\cite{kn:PeS1} \cite{kn:PeS2}, \cite{kn:PeS3}).

In what follows we consider the following {\it lower unstable pinching condition} for
$\phi_t$ and $\mt$:

\ms

\noindent
{\sc (LUPC)}:  {\it There exist  constants $C > 0$ and $0 < \alpha \leq \beta < \alpha_2 \leq \beta_2$,
and for  every $x\in \mt$ constants $\alpha_1(x) \leq \beta_1(x)$ with $\alpha \leq \alpha_1(x) 
\leq \beta_1(x) \leq \beta$
and $2\alpha_1(x) - \beta_1(x) \geq \alpha$ and a $d\phi_t$-invariant splitting 
$E^u(x) = E^u_1(x)\oplus E^u_2(x)$,
continuous with respect to $x \in \mt$, such that
$$\frac{1}{C} \, e^{\alpha_1(x) \,t}\, \|u\| \leq \| d\phi_{t}(x)\cdot u\| \leq C\, e^{\beta_1(x)\,t}\, \|u\|
\quad, \quad  u\in E^u_1(x) \:\:, t > 0 \;,$$
and}
$$\frac{1}{C} \, e^{\alpha_2 \,t}\, \|u\| \leq \| d\phi_{t}(x)\cdot u\| \leq C\, e^{\beta_2\,t}\, \|u\|
\quad, \quad  u\in E^u_2(x) \:\:, t > 0 \;.$$

\ms

In (LUPC) the lower part of the spectrum of $d\phi_t$ over $E^u$ is (point-wisely) pinched, however 
there is no restriction on the rest of the spectrum, except that it should be uniformly separated 
from the lower part.

Under the above condition the distribution $E^u_2(x)$ ($x\in \mt$) is integrable 
(see e.g. \cite{kn:Pes}), so (assuming $\ep_0 > 0$ is
small enough) there exists a $\phi_t$-invariant family $W^{u,2}_{\ep_0}(x)$ ($x\in \mt)$ of $C^2$ 
submanifolds of $W^u_{\ep_0}(x)$ such that $T_x(W^{u,2}_{\ep_0}(x)) = E^u_2(x)$ for all $x\in \mt$. 
Moreover (see Theorem 6.1 in \cite{kn:HPS} or the proof of Theorem B in \cite{kn:PSW}), for any 
$x\in \mt$, the map 
$\mt \cap W^u_{\ep_0}(x) \ni y \mapsto E^u_2(y)$ is $C^1$. However in general the distribution  
$E^u_1(x)$ ($x\in \mt$)  does not have to be integrable (see  \cite{kn:Pes}).

We now make the  additional assumption that  $E^u_1(x)$ ($x\in \mt$) is integrable:

\ms

\noindent
{\sc (I)}:  {\it There exist $\ep_0 > 0$ and a continuous $\phi_t$-invariant family $W^{u,1}_{\ep_0}(x)$ 
($x\in \mt)$ of $C^2$ submanifolds
of $W^u_{\ep_0}(x)$ such that $T_x(W^{u,1}_{\ep_0}(x)) = E^u_1(x)$ for all $x\in \mt$, and moreover 
for any $\ep > 0$ and 
any $x\in \mt$, $\mt\cap W^u_{\ep}(x)$ is {\bf not} contained in  $W^{u,2}_{\ep}(x)$. }

\ms


Our main result in this paper is the following.

\bs

\noindent
{\sc Theorem 1.1.} {\it Let $\phi_t$ and $\mt$ satisfy the conditions {\rm (LUPC)} 
and {\rm (I)}. Then $\phi_t$
has regular distortion along unstable manifolds over $\mt$.}

\bs

A simplified case of the above is presented by the following {\it pinching condition}:

\ms

\noindent
{\sc (P)}:  {\it There exist  constants $C > 0$ and $\beta \geq \alpha > 0$ such that for every 
$x\in \mt$ we have
$$\frac{1}{C} \, e^{\alpha_x \,t}\, \|u\| \leq \| d\phi_{t}(x)\cdot u\| \leq C\, e^{\beta_x\,t}\, \|u\|
\quad, \quad  u\in E^u(x) \:\:, t > 0 \;,$$
for some constants $\alpha_x, \beta_x > 0$ depending on $x$ but independent of $u$ and $t$ with
$\alpha \leq \alpha_x \leq \beta_x \leq \beta$ and $2\alpha_x - \beta_x \geq \alpha$ for all $x\in \mt$.}

\ms

Clearly the condition (P) is (LUPC) in the special case when $E^u_2(x) = 0$ for all $x\in \mt$.  
Notice that when the local unstable manifolds are one-dimensional the condition (P) is always satisfied. 
In higher dimensions a well-known example when (P) holds is the geodesic flow on a manifold with 
strictly negative sectional curvature satisfying the so called $\frac{1}{4}$-pinching condition 
(see \cite{kn:HP}). For open billiard flows (in any dimension) it was shown in \cite{kn:St2} that 
if the distance between the scatterers is large compared with the maximal sectional curvature of 
the boundaries, then the condition (P) is satisfied over the non-wandering set.

As a special case of Theorem 1.1, we  get the following.

\bs

\noindent
{\sc Proposition 1.2.} {\it Let $\phi_t$ and $\mt$ satisfy the conditions {\rm (P)}.
Then $\phi_t$ has regular distortion along unstable manifolds over $\mt$.}

\bs

Our strategy in this paper is to prove Proposition 1.2 first. This is done in sections 3 and 4 below.
Then in section 5 we generalize the arguments from sections 3 and 4 to prove Theorem 1.1. While
the latter is technically more difficult, the main ideas in its proof are almost the same as those
in the proof of Proposition 1.2.

In section 3  we prove that under the condition (P) the flow $\phi_t$ on unstable manifolds is 
locally conjugate  to its linearization $d\phi_t$  via $C^1$ (local) maps with uniformly bounded 
derivatives.  Linearization via a family of homeomorphisms exists in 
the general case (see section 4 in \cite{kn:PS}), however it is not clear whether one can make 
it Lipschitz without any additional conditions. The arguments used to prove some recent results 
on smoothness of the linearizing homeomorphism in Hartman-Grobman type theorems do not seem to 
be easily applicable  to vector bundles (see \cite{kn:GHR}  and the references there). 
The linearization from section 3 is used in section 4 to derive Proposition 1.2. 

In section 5 we use slight modifications of the arguments from sections 3 and 4 
to prove Theorem 1.1. It should be stressed that the central  
part of the arguments in sections 4 and 5 is to establish a local version of regular distortion along 
unstable manifolds, where  e.g. (1.1) is satisfied at a single point $z\in \mt$ (with a constant
$R = R(z,\delta, \ep)$ depending on $z$ as well) -- see Lemma 4.2 below. It is not difficult to 
see (using a variation of the arguments  in sections 3-5) that a similar local result can be proved 
at Lyapunov regular points $z\in \mt$ for any $C^2$ flow  on any basic set\footnote{Now the role of the 
'bottom of the unstable spectrum' is played by the exponential of the least positive Lyapunov exponent.}. 
However, it is not clear how one can get a uniform global result over $\mt$ from such local results.

\section{Preliminaries}
\setcounter{equation}{0}

Throughout this paper $M$ denotes a $C^2$ complete (not necessarily compact) 
Riemann manifold,  and $\phi_t : M \longrightarrow M$ ($t\in \R$) a $C^2$ flow on $M$. A
$\phi_t$-invariant closed subset $\mt$ of $M$ is called {\it hyperbolic} if $\mt$ contains
no fixed points  and there exist  constants $C > 0$ and $0 < \lambda < 1$ such that 
 there exists a $d\phi_t$-invariant decomposition 
$T_xM = E^0(x) \oplus E^u(x) \oplus E^s(x)$ of $T_xM$ ($x \in \mt$) into a direct sum of non-zero linear subspaces,
where $E^0(x)$ is the one-dimensional subspace determined by the direction of the flow
at $x$, $\| d\phi_t(u)\| \leq C\, \lambda^t\, \|u\|$ for all  $u\in E^s(x)$ and $t\geq 0$, and
$\| d\phi_t(u)\| \leq C\, \lambda^{-t}\, \|u\|$ for all $u\in E^u(x)$ and  $t\leq 0$.

The flow $\phi_t$ is called an {\it Axiom A flow} on $M$ if the non-wandering set 
of $\phi_t$  is a disjoint union of a finite set  consisting of  fixed
hyperbolic points and a compact hyperbolic subset containing no fixed points in which
the periodic points are dense (see e.g. \cite{kn:KH}). 
A non-empty compact $\phi_t$-invariant hyperbolic subset $\mt$ of $M$ which is not a single 
closed orbit is called a {\it basic set} for $\phi_t$ if $\phi_t$ is transitive on $\mt$ 
and $\mt$ is locally maximal, i.e. there exists an open neighbourhood $V$ of
$\mt$ in $M$ such that $\mt = \cap_{t\in \R} \phi_t(V)$. When $M$ is compact and $M$ itself
is a basic set, $\phi_t$ is called an {\it Anosov flow}.

For $x\in \Lambda$ and a sufficiently small $\epsilon > 0$ let 
$$\wloc^s(x) = \{ y\in M : d (\phi_t(x),\phi_t(y)) \leq \epsilon \: \mbox{\rm for all }
\: t \geq 0 \; , \: d (\phi_t(x),\phi_t(y)) \to_{t\to \infty} 0\: \}\; ,$$
$$\wloc^u(x) = \{ y\in M : d (\phi_t(x),\phi_t(y)) \leq \epsilon \: \mbox{\rm for all }
\: t \leq 0 \; , \: d (\phi_t(x),\phi_t(y)) \to_{t\to -\infty} 0\: \}$$
be the (strong) {\it stable} and {\it unstable manifolds} of size $\epsilon$. Then
$E^u(x) = T_x \wloc^u(x)$ and $E^s(x) = T_x \wloc^s(x)$. 
Given $\delta > 0$, set $E^u(x;\delta) = \{ u\in E^u(x) : \|u\| \leq \delta\}$;
$E^s(x;\delta)$ is defined similarly. For any $A\subset M$ and  $I \subset \R$  denote
$\phi_I(A) = \{\; \phi_t(y)\; : \; y\in A, t \in I \; \}.$

It follows from the hyperbolicity of $\mt$  that if  $\epsilon_0 > 0$ is sufficiently small,
there exists $\ep_1 > 0$ such that if $x,y\in \mt$ and $d (x,y) < \ep_1$, 
then $W^s_{\ep_0}(x)$ and $\phi_{[-\ep_0,\ep_0]}(W^u_{\ep_0}(y))$ intersect at exactly 
one point $[x,y ] \in \mt$  (cf. \cite{kn:KH}). That is, there exists a unique $t\in [-\ep_0, \ep_0]$ such that
$\phi_t([x,y]) \in W^u_{\ep_0}(y)$. Setting $\Delta(x,y) = t$, defines the so called {\it temporal distance
function}. For $x, y\in \mt$ with $d (x,y) < \ep_1$, define
$\pi_y(x) = [x,y] = W^s_{\ep}(x) \cap \phi_{[-\ep_0,\ep_0]} (W^u_{\ep_0}(y))\;.$
Thus, for a fixed $y \in \mt$, $\pi_y : W \longrightarrow \phi_{[-\ep_0,\ep_0]} (W^u_{\ep_0}(y))$ is the
{\it projection} along local stable manifolds defined on a small open neighbourhood $W$ of $y$ in $\mt$.
Choosing $\ep_1 \in (0,\ep_0)$ sufficiently small, 
the restriction
$\pi_y: \phi_{[-\ep_1,\ep_1]} (W^u_{\ep_1}(x)) \longrightarrow \phi_{[-\ep_0,\ep_0]} (W^u_{\ep_0}(y))$
is called a local stable holonomy map\footnote{In a similar way one can define
holonomy maps between any two sufficiently close local transversals to stable laminations; see e.g.
\cite{kn:PSW}.}. Combining it with a shift along the flow we get another {\it local stable holonomy  map}
$\hhs_{x,y} : W^u_{\ep_1}(x) \cap \mt  \longrightarrow W^u_{\ep_0}(y) \cap \mt$.
In a similar way one defines local holonomy maps along unstable laminations.


\def\hcc{\widehat{\cc}}
\def\Bmt{\overline{B_{\ep_0}(\mt)}}
\def\Lye{L_{y,\eta}}
\def\Lyep{L^{(p)}_{y,\eta}}
\def\Fyp{F^{(p)}_y}
\def\Fxp{F^{(p)}_x}
\def\Lxx{L_{x,\xi}}
\def\Lxxp{L^{(p)}_{x,\xi}}

\section{Linearization of pinched Axiom A flows}
\setcounter{equation}{0}

Let $M$ be a $C^2$ complete Riemann manifold, $\phi_t$  be a  $C^2$  flow on $M$, and let 
$\mt$ a basic set for $\phi_t$. In this section we assume that $\phi_t$ and $\mt$ satisfy the 
pinching condition (P) from the Introduction.

We  prove that under this condition  the flow $\phi_t$ on unstable manifolds is locally conjugate 
to its  linearization  $d\phi_t$  via $C^1$ (local) maps with uniformly bounded derivatives.
This is used in section 4 to show that
the condition (P) implies regular distortion along unstable manifolds over $\mt$.



Fix constants $C > 0$ and $\beta \geq \alpha > 0$ and for each $x\in \mt$ constants $\alpha_x \leq \beta_x$ 
with the properties in (P). Throughout we use the notation $m(A) = 1/\|A^{-1}\|$ for an invertible 
linear operator $A$. 

Take $\ep_0 \in (0,1/2)$  such that for any $x\in \mt$,
$\exp^u_x : E^u (x; \ep_0) \longrightarrow \exp^u_x (E^u (x; \ep_0)) \subset W^u_{\ep_0}(x)$
is a diffeomorphism. We will assume the constant $C \geq 1$ is so large that
\be
\|d \exp^u_x(u)\| \leq C \quad ,\quad  \|(d\exp^u_x(u))^{-1}\| \leq C \quad, 
\quad x\in W^u_{\ep_0}(\mt) \:, \: u \in E^u (x ; \ep_0)\;.
\ee

Choose an arbitrary $t_0 > 0$ such that
\be
\gamma = 8 C^3\, e^{-\alpha\, t_0/2} < 1\;,
\ee
and fix it. Set $f = \phi_{t_0}$ and fix a constant $\ep_1 \in (0,\ep_0]$ such that for any 
$x\in W^u_{\ep_0}(\mt)$ the map
$$\hf_x = (\exp^u_{f(x)})^{-1} \circ f\circ \exp^u_x : E^u (x ; \ep_1) \longrightarrow E^u (f(x); \ep_0)\;$$
is well-defined (and therefore $C^2$).

For any $y \in W^u_{\ep_0}(\mt)$ and any integer $k \geq 1$ we will use the notation
$$\hf_y^k = \hf_{f^{k-1}(y)} \circ \ldots \circ \hf_{f(y)} \circ \hf_y\quad,
\quad \hf_y^{-k} = (\hf_{f^{-k}(y)})^{-1} \circ \ldots \circ (\hf_{f^{-2}(y)})^{-1} \circ 
(\hf_{f^{-1}(y)})^{-1} \;,$$
at any point where these sequences of maps are well-defined.

Given $\ep \in (0,\ep_1)$, the set
$$\hLa_x(\ep) = \{ u\in E^u(x;\ep) : \exp^u_x(u) \in \mt\}$$
is the {\it local representative of $\mt$} in $E^u(x)$.  Notice that   
$\hf^{-1}_x (\hLa_x) \subset \hLa_{f^{-1}(x)}$ for any $x\in \mt$.

Our aim in this section is to show that under the pinching condition (P) we can locally linearize the maps
$\hf_x$ by a continuous family of local diffeomorphisms and that family `linearizes' the sets $\hLa_x$, 
as well. More precisely, we have the following:

\bs

\noindent
{\sc Theorem 3.1.} {\it Assume that $\phi_t$ and $\mt$  satisfy the pinching condition {\rm (P)}. Then
there exists a constant $\ep_2 \in (0,\ep_1/2]$ such that for every $x \in W^u_{\ep_2}(\mt)$ we have
the following:}

\ms

(a) {\it For  every $u\in E^u(x; \ep_2)$ there exists
$\di F_x(u) = \lim_{p\to\infty} d\hf^p_{f^{-p}(x)}(0)\cdot \hf_x^{-p}(u) \in E^u(x; 2\ep_2)\;.$
Moreover,  there exists a constant $C_1 > 0$ such that
$\|F_x(u) - d\hf^p_{f^{-p}(x)}(0)\cdot \hf_x^{-p}(u)\| \leq C_1\, \gamma^p\, \|u\|^2$ 
for any $u\in E^u(x,\ep_2)$ and any integer $p \geq 0$.}

\ms 

(b) {\it The map $F_x : E^u( x;\ep_2) \longrightarrow F_x (E^u( x;\ep_2)) \subset E^u(x; 2\ep_2)$ 
is a $C^1$ diffeomorphism with uniformly bounded derivatives. }

\ms

(c) {\it For any  integer $q\geq 1$  we have 
$d\hf_{x}^{-q}(0) \circ F_{x} (v) = F_{ f^{-q}(x)} \circ \hf_{x}^{-q} (v)$ for any $v \in E^u (x;\ep_2)$.}

\ms

(d) {\it For  any $\xi, u \in E^u(x;\ep_2/2)$ there exist the limits
$$\di \Lxx  = \lim_{p\to\infty} d\hf^p_{f^{-p}(x)}(\hf_x^{-p}(\xi))\circ d\hf_x^{-p}(0) \quad , \quad 
F_{x,\xi}(u) = \lim_{p\to\infty} d\hf^p_{f^{-p}(x)}(\hf_x^{-p}(\xi))\circ \hf_x^{-p}(u)\;.$$
Moreover, for the linear map $\Lxx$ we have  $\|\Lxx\| \leq 2$, while $F_{x,\xi}(u) = \Lxx \circ F_x (u)$. }

\bs

As an immediate consequence of the above one gets the following, where we use the notation
$ \hphi_{x,t} = (\exp^u_{\phi_t(x)})^{-1} \circ \phi_t\circ \exp^u_x$ for any $x\in K$ and $t\in \R$.

\bs

\noindent
{\sc Corollary  3.2.} {\it  For  any $x\in W^u_{\ep_2}(\mt)$ we have 
$F_x(u) = \lim_{t\to \infty} d\phi_t(\phi_{-t}(x))\cdot \hphi_{x,t}(u)$. Moreover,
$d\phi_t(x) \cdot F_x (u) = F_{\phi_t(x)} (  \hphi_{x,t} (u))$ for $t \geq 0$ and $u\in E^u(x)$ with $\| \hphi_{x,t} (u)\| \leq \ep_2$. }

\bs

The rest of this section is devoted to the proof of Theorem 3.1.

Taylor's formula and the compactness of $W^u_{\ep_0}(\mt)$,  imply that there exists a constant 
$D > 0$ such that 
\be
\|\hf_x(v) - \hf_x(u) - d\hf_x(u) \cdot (v-u)\| \leq D\, \|v-u\|^{2} \quad, 
\quad x\in W^u_{\ep_0}(\mt)\:, \: u,v \in E^u (x ; \ep_1)\;.
\ee
Since $\phi_t$ is $C^2$,  we can take $D$ so large that 
\be
\|d\hf_x(u) - d\hf_x(0)\| \leq D\, \|u\| \quad , \quad x\in W^u_{\ep_0}(\mt) \:, \: u \in E^u (x;\ep_1)\;.
\ee
Combining the latter with (3.3) gives
\be
\|\hf_x(v) - \hf_x(u) - d\hf_x(0) \cdot (v-u)\| \leq D\, \left[\|v-u\|^{2} + \|u\|\, \|v-u\|\right]
\ee
for all $x\in W^u_{\ep_0}(\mt)$ and all $u, v \in E^u (x ; \ep_1)$.

Next, for  each $x\in \mt$ set  $\hal_x = \alpha_x - \alpha/8$ and 
$\hbe_x =  \beta_x + \alpha/4$. Then $2\hal_x -  \hbe_x  \geq \alpha/2$, so by (3.2),
$$8C^3\, e^{(\hbe_x- 2 \hal_x)\, t_0} \leq  \gamma = 8C^3\, e^{-\alpha\, t_0/2} < 1\;.$$

Assuming  $\ep_1 <  \frac{C\, e^{\alpha \, t_0}}{D} \, 
\min\{ e^{\alpha\, t_0/16} -1 , (1 - e^{-\alpha\, t_0/8})/ C^2\}\;,$
one derives from (P) and (3.4) that 
\be
\frac{1}{C}\, e^{(\hal_x + \alpha/16)\, t_0}\, \|u\| \leq \| d\hf_x(\eta)\cdot u\| 
\leq C\, e^{(\hbe_x - \alpha/8)\, t_0}\, \|u\| 
\ee
for all $x\in \mt$, $\eta\in E^u(x;\ep_1)$ and  $u\in E^u(x)$.

Finally, assuming also that $\ep_1 \leq \frac{1}{4 D\, C}$, the above and (3.3) imply
$$\mu_x\,e^{\alpha\, t_0/16}\, \|u-v\| \leq \|\hf_x(v) - \hf_x(u)\| \leq 
\lambda_x\,e^{-\alpha\, t_0/8}\,  \|u-v\| 
\quad, \quad x\in \mt\:, \: u,v \in E^u (x ; \ep_1)\;,$$
where 
$1 < \mu_x = \frac{1}{2C}\, e^{\hal_x \,t_0} < \lambda_x =  2C\, e^{\hbe_x \,t_0}\;.$

In fact, replacing $\ep_1 > 0$ and $\ep_2 > 0$ by smaller numbers if necessary, we can arrange that
\be
\mu_x\, \|u-v\| \leq \|\hf_x(v) - \hf_x(u)\| \leq \lambda_x\, \|u-v\| 
\quad, \quad x\in W^u_{\ep_2}(\mt);  u,v \in E^u (x ; \ep_1)\;.
\ee
Indeed, assume $\ep_1 > 0$ and $\ep_2 > 0$ are so small that
$G_x^y = (\exp^u_y)^{-1}\circ \exp^u_x : E^u(x;\ep_1) \longrightarrow E^u(y;\ep_0)$
is well-defined for $y\in W^u_{\ep_2}(x)$, $x \in \mt$; it is then a $C^2$ map with uniformly 
bounded derivatives. 
Since $G_x^x = \id$, $dG_x^y$ can be made arbitrarily close to $\id$ taking $\ep_2 > 0$ sufficiently small.
Fix $\delta > 0$ so small that $(1+\delta)^2 < e^{\alpha\, t_0/8}$ and 
$(1-\delta)^2 > e^{-\alpha\, t_0/16}$, and then take $\ep_1' > 0$ and $\ep_2 > 0$ so small that
$\|dG_x^y - \id\| \leq \delta$ and $\|G_x^y(u)\| \leq \ep_1$ for all $x\in \mt$,  
$y\in W^u_{\ep_2}(x)$ and $u\in E^u(x;\ep_1')$. 
Given  $x_0\in \mt$,  $x\in W^u_{\ep_2}(x_0)$ and  $u,v \in E^u(x;\ep'_1)$, setting 
$u' = G_x^{x_0}(u)$, $v' = G_x^{x_0}(v)$, 
it is easy to see that 
$\|\hf_x(u)- \hf_x(v)\| 
\leq \lambda_x\, \|u-v\|$.
Similarly, $\|\hf_x(u)- \hf_x(v)\| \geq \mu_x\, \|u-v\|$. Now replacing $\ep_1$ by $\ep_1'$ proves (3.7).

Notice that
\be
 \lambda_x \, \mu_x^{-2} = 8 C^3\, e^{(\hbe_x - 2\hal_x)t_0} \leq \gamma  < 1\quad , \quad x\in \mt\;.
\ee

The following is the main step in the proof of Theorem 3.1.

\bs

\noindent
{\sc Lemma 3.3.} {\it There exist constants $C_1 > 0$, $\ep_1 \in (0,\ep_0]$ and 
$\ep_2 \in (0,\ep_1/(2C_1)]$ with the following properties:}

\ms

(a) {\it If $z\in W^u_{\ep_2}(\mt)$ and $\|\hf^p_z(v)\| \leq \ep_2$ for some $v\in E^u (z ; \ep_1)$ 
and some integer 
$p \geq 1$, then $\|d\hf_z^p(0)\cdot v - \hf^p_z(v)\| \leq C_1\, \|\hf^p_z(v)\|^2$, and therefore
$\|d\hf^p_z(0)\cdot v\| \leq 2 \|\hf^p_z(v)\|$.
Similarly, if  $\|d\hf^p_z(0)\cdot v\|\leq \ep_2$ for some $v\in E^u (z; \ep_1)$ and 
some integer $p \geq 1$, then $\|\hf_z^p(v) - d\hf^p_z (0)\cdot v\| \leq C_1\, \|d\hf^p_z(0)\cdot v\|^2$, 
and so $\|\hf^p_z(v) \| \leq 2 \|d\hf^p_z(0)\cdot v\|$.}

\ms

(b)  {\it If $z\in W^u_{\ep_2}(\mt)$ and $\|\hf^p_z(v)\| \leq \ep_2$, $\|\hf^p_z(\xi)\| \leq \ep_2$ for some 
$v, \xi\in E^u (z ; \ep_2)$ and some integer $p \geq 1$, then 
$\|d\hf_z^p(\xi)\cdot v - d\hf^p_z(0) \cdot v\| \leq C_1\, \|d\hf^p_z(0)\cdot v\|\, \|\hf_z^p(\xi)\|$.}

\ms

(c)  {\it For any $y \in W^u_{\ep_2}(\mt)$ and any integer $p \geq 1$ the map
$\Fyp = d \hf^p_{f^{-p}(y)}(0) \circ (\hf^p_y)^{-1}:  E^u_y (\ep_2) \longrightarrow E^u_y (2\ep_2)$  
is such that
$$\left\| \left[ \Fyp (a) - \Fyp (b) \right] - [a - b ]\right\|  
\leq C_1 \, \left[ \|a - b \|^2 + \| b \|\cdot \|a - b\| \right] 
\quad, \quad a,b \in E^u (y ; \ep_2)\;.$$
Similar estimates hold for the map
$\hf^p_{f^{-p}(y)} \circ (d\hf^p_y (0))^{-1}:  E^u (y ; \ep_2 ) \longrightarrow E^u (y ; 2\ep_2)$.}

\ms

(d) {\it For any $y \in W^u_{\ep_2}(\mt)$, any $\eta \in E^u(y;\ep_2)$ and any integer 
$p \geq 1$ the linear map\\
$L^{(p)}_{y,\eta} = d \hf^p_{f^{-p}(y)} (\hf_x^{-p}(\eta)) \circ (d\hf^{-p}_y(0)):  
E^u (y) \longrightarrow E^u (y)$  is such that
$\left\| \Lyep (a)  - a \right\|  \leq C_1 \,  \|a  \|\, \| \eta \|$ for all $a \in E^u (y ; \ep_2)$,
so $\|\Lyep\| \leq 2$ for all $p$.}

\bs

\noindent
{\it Proof of Lemma} 3.3.  Set $C_1 = \frac{10\,  D }{1-\gamma}$ and choose $\ep_2 \in (0,\ep_1/(2C_1)]$ 
with (3.7) and
$\frac{30\,  D\, \ep_2 }{1-\gamma} < \frac{1}{2}$.

\ms

(a) Fix arbitrary $z_0\in \mt$ and $z\in W^u_{\ep_0}(z_0)$, and let  $v \in E^u (z;\ep_1)$ and 
$p \geq 1$ be such that
$\|\hf^p_z(v)\| \leq \ep_2$. Set  $z_j = f^j(z)$, $v_j = \hf^j_z(v) \in E^u (z_j)$ and 
$w_j = d\hf_z^j(0)\cdot v \in E^u (z_j)$.

Then (3.7)  implies
\be
\|\hf^k_z(v)\| \leq \frac{1}{\mu_{z_k}}\, \|\hf_{z_k}(\hf_z^k(v))\| 
=  \frac{1}{\mu_{z_k}}\, \|\hf_z^{k+1}(v)\|\leq \ldots \leq 
 \frac{1}{\mu_{z_k}\, \mu_{z_{k+1}} \ldots \mu_{z_{p-1}} }\, \|\hf_z^p(v)\|\;
 \ee
 for all $ k = 0,1, \ldots,p-1$.

By (3.3), $\|\hf_z(v) - d\hf_z(0)\cdot v\| \leq D \, \|v\|^{2}$, so
$w_1 = d\hf_z(0)\cdot v = \hf_z (v) + u_1$
for some $u_1 \in E^u (z_1)$ with $\|u_1\|\leq D\, \|v\|^{2}$. Hence
\be
d\hf_z^2(0)\cdot v = d\hf_{z_1}(0)\circ d\hf_z (0)\cdot v
= d\hf_{z_1}(0) \cdot (\hf_z (v)) + d\hf_{z_1}(0)\cdot u_1\;.
\ee
Using (3.3) again, 
$\|\hf_{z_1}((\hf_z (v))) - d\hf_{z_1}(0)\cdot \hf_z (v)\| \leq D \, \|\hf_z (v)\|^{2}\;,$
so $d\hf_{z_1}(0)\cdot (\hf_z (v)) = v_2 + u_2$
for some $u_2\in E^u (z_2)$ with $\|u_2\| \leq D \, \|v_1\|^{2}$. Now (3.10) gives
$w_2 = v_2 + u_2 + d\hf_{z_1}(0)\cdot u_1\;.$

In this way one proves by induction that for any  $k = 1, \ldots, p$ we have
\be
w_k = v_k + u_k + d\hf_{z_{k-1}}(0)\cdot u_{k-1} + d\hf^2_{z_{k-2}}(0)\cdot u_{k-2} +
\ldots + d\hf_{z_1}^{k-1}(0)\cdot u_1\;,
\ee
where $u_j \in E^u (z_j)$ and $\|u_j\| \leq D \, \|v_{j-1}\|^{2}$ for all $j = 1,\ldots,k-1$.

Next, (3.9) implies
$\|u_j\| \leq D\, \|v_{j} \|^{2} \leq  \frac{D}{\mu^2_{z_{j}}\, \mu^2_{z_{j+1}} 
\ldots \mu^2_{z_{p-1}} }\,  \|\hf^p_z(v)\|^2\;.$
Combining the latter with (3.6) and (3.8) gives
$\|d\hf_{z_j}^{p-j}(0)\cdot u_j\|  \leq  \lambda_{z_j} \, \lambda_{z_{j+1}} \ldots 
\lambda_{z_{p-1}} \, \|u_j\|  \leq   D \, \gamma^{p-j} \, \|\hf^p_z (v)\|^2$
for all $j = 1, \ldots, p$. It then follows from (3.11) with $k = p$  that
$$\|d\hf_z^p(0)\cdot v - \hf^p_z (v) \| = \|w_p - v_p\|  \leq  D \, \| \hf^p_z (v)\|^2 \, 
\sum_{j=1}^p \gamma^{p-j} \leq C_1 \, \| \hf^p_z (v)\|^2 \;.$$
According to the choice of $\ep_2$, the latter implies $\|d\hf_z^p(0)\cdot v\| \leq  2\| \hf^p_z (v)\|$. 

To prove the second half of part (a),  assume that $\|d\hf_z^p(0)\cdot v\| \leq \ep_2$. Then (3.6) implies
$ \|d\hf_z^p(0)\cdot v\| = \|d\hf^{p-j}_{z_j}(0)\cdot w_j\| \geq  
\mu_{z_{j}}\, \mu_{z_{j+1}} \ldots \mu_{z_{p-1}} \, \|w_j\|\;,$
so 
\be
\|w_j\| \leq  \frac{1}{\mu_{z_{j}}\, \mu_{z_{j+1}} \ldots \mu_{z_{p-1}} } \, \|d\hf_z^p(0)\cdot v\| 
 \quad, \quad 0 \leq j \leq p-1\;.
\ee

We will show by induction on $k$ that $v_k = \hf^k_z (v)$ is well-defined and
$v_k \in E^u (z_k ;\ep_1)$ for all $k = 0,1, \ldots, p$.  It follows from  (3.3) and (3.4) that
$$\|w_{k+1} - \hf_{z_k}(w_k)\| = \| d\hf_{z_k}(0)\cdot w_k - \hf_{z_k}(w_k)\| \leq D \, \|w_k\|^{2}
\leq \frac{D}{ \mu^2_{z_{k}}\, \mu^2_{z_{k+1}} \ldots \mu^2_{z_{p-1}} }\, \|d\hf_z^p(0)\cdot v\|^2\;.$$
This and (3.8) yield (showing in the meantime by induction that $\hf_{z_k}^{p-k}(w_k)\in E^u (z_p; \ep_1)$ 
for all $k$)
\begin{eqnarray*}
\|\hf_{z_{k+1}}^{p-k-1}(w_{k+1}) - \hf_{z_k}^{p-k}(w_k)\|
& \leq &  \lambda_{z_{k+1}}\, \lambda_{z_{k+2}} \ldots \lambda_{z_{p-1}}\,\|w_{k+1} - \hf_{x_k}(w_k)\| \\
& \leq & D\, \frac{\lambda_{z_k}}{\mu^2_{z_k}}\,  \frac{\lambda_{z_{k+1}}}{\mu^2_{z_{k+1}}}
\ldots  \frac{\lambda_{z_{p-1}}}{\mu^2_{z_{p-1}}}\, \|d\hf_z^p(0)\cdot v\|^2 \leq
 D\, \gamma^{p-k}\, \|w_p\|^2\;.
\end{eqnarray*}
Hence
$ \|w_p - v_p \| \leq \sum_{k=0}^{p-1} \|\hf_{z_{k+1}}^{p-k-1}(w_{k+1}) - \hf_{z_k}^{p-k}(w_k)\|
\leq D\, \|w_p\|^2 \, \sum_{k=0}^{p-1}  \gamma^{p-k} \leq C_1\, \|w_p\|^2  \;,$
provided $0 < C_1 \leq \frac{D}{1-\gamma}$. The above also implies $\|v_p\| \leq 2\|w_p\|$.

\ms

The proofs of parts (b) and (c) are essentially repetitions of the proof of part (a), so we omit the details. 

\ms

(d) Given $a \in E^u(y,\ep_2)$, set $v = d\hf_y^{-p}(0)\cdot a$. Using part (b) with $z = f^{-p}(y)$ and $\xi = \hf^{-p}(\eta)$,
we get
$$\|L^{(p)}_{y,\eta}(a) - a\| = \| d\hf_z^p(\xi) \cdot v - d \hf_z^p(0)\cdot v\| \leq C_1\, \|a\|\, \|\eta\|\,$$
which proves the claim.
\endofproof

\bs

\noindent
{\it Proof of Theorem 3.1.} Let $\ep_2 \in (0,\ep_1/(2C_1)]$ and $C_1 > 1$ be as in Lemma 3.3. Fix
arbitrary $x_0\in \mt$ and $x\in W^u_{\ep_2}(x_0)$  and set $x_p = f^{-p}(x)$ for any integer $p \geq 0$.
(Notice the different meaning of this notation here.)
In what follows we use the maps $\Fyp$ ($y\in \mt$, $p \geq 1$) from  Lemma 3.3.

\ms

(a), (d) Given $u, \xi \in E^u(x, \ep_2/2)$ and $p \geq 1$, set $u_p = \Fxp(u) \in E^u(x;2\ep_2)$, 
$\zeta_p = \Lxxp (u)$
To show that the sequences $\{ u_p\}$ and $\zeta_p$ are  Cauchy, consider any $q > p$ and set
$v = \hf_x^{-p}(u) \in E^u(x_p, \ep_2)$. By (3.7),
$\|v\| \leq \frac{\|u\|}{\mu_{x_{1}} \, \mu_{x_2} \ldots \mu_{x_{p}} }\,\;.$
From Lemma 3.3(a) we know that  $\|v_{q-p} - v\| \leq C_1\, \|v\|^2$, i.e.
$\|d \hf_{x_q}^{(q-p)}(0)\cdot (\hf^{-(q-p)}_{x_p}(v)) - v\| \leq C_1\, \|v\|^2\;.$
Applying $d\hf_{x_p}^p(0)$ to the latter and using the estimate for $\|v\|$ and (3.8), we get
$$\|u_q - u_p\| = \|d\hf^q_{x_q}(0)\cdot (\hf^{-q}_x(u)) - d\hf^p_{x_p}(0)\cdot v\|
\leq \lambda_{x_p} \, \lambda_{x_{p-1}}\ldots \lambda_{x_1}\, C_1\, \|v\|^2 \leq C_1\, \gamma^p\, \|u\|^2\;.$$
Thus, $\{ u_p\}$ is Cauchy, so there exists $F_x(u) = \lim_{p\to\infty} u_p$. Moreover, letting 
$q\to \infty$ in the 
above gives $\|F_x(u) - u_p\| \leq  C_1\, \gamma^p\, \|u\|^2$ for all  $u\in E^u(x;\ep_2)$ and $p\geq 1$.

In a similar way, $\|\zeta_q - \zeta_p\| \leq  C_1\, \gamma^p\, \|u\| \, \|\xi\|$, so $\{ \zeta_p\}$ 
is Cauchy. Thus, 
there exists $\Lxx(u) = \lim_{p\to\infty} \Lxxp (u)$. Moreover, letting $q\to \infty$ in the 
above gives $\|\Lxx(u) - \Lxxp (u)\| \leq  C_1\, \gamma^p\, \|u\|\, \|\xi\|$ for all  
$u\in E^u(x;\ep_2/2)$ and $p\geq 1$.
By Lemma 3.3(d), $\|\Lxxp\|\leq 2$ for all $p \geq 1$, so $\|\Lxx\| \leq 2$, as well.

It remains to show that  $F_{x,\xi}(u)$ exists and $F_{x,\xi}(u) = \Lxx \circ F_x (u)$. 
Setting $\xi_p = \hf^{-p}_x(\xi)$, we have 
$$d\hf^p_{x_p}(\xi_p)\cdot \hf^{-p}_x(u) = \Lxxp \circ \left(d\hf^p_{x_p}(0)\cdot \hf^{-p}_x(u)  \right)
= \Lxxp \circ \Fxp (u)\;. $$
Next,
\begin{eqnarray*}
\|\Lxx(F_x(u)) - \Lxxp (\Fxp(u))\|
& \leq & \|\Lxx(F_x(u)) - \Lxxp (F_x(u))\| + \|\Lxxp(F_x(u)) - \Fxp(u))\|\\
& \leq & C_1\, \gamma^p\,\| F_x(u)\|\, \|\xi\| + 2 \|F_x(u) - \Fxp(u)\| \to 0
\end{eqnarray*}
as $p \to \infty$. Thus, there exists $\lim_{p\to \infty} \Lxxp(\Fxp(u)) = \Lxx(F_x(u))$. 

\ms

(b) Given $u, v\in E^u(x;\ep_2)$ and $p \geq 0$, it follows from Lemma 3.3(c) that
$$\|(u_p - v_p) - (u-v)\| \leq C_1\, [ \|u-v\|^2 + \|v\|\cdot \|u-v\|]\;.$$
Letting $p \to \infty$, gives $\|F_x(u) - F_x(v) - (u-v)\| \leq C_1\, [\|u-v\|^2 + \|v\|\cdot \|u-v\|] $. 
In particular, there exists $dF_x(0) = I = \id$.

It remains to show that $F_x$ is a $C^1$ diffeomorphism.
Assuming $\ep_1 \in (0,\ep_0]$ and $\ep_2 \in (0,\ep_1/(2C_1)]$ are small enough, 
$G_x^y = (\exp^u_y)^{-1}\circ \exp^u_x : E^u(x;\ep_1) \longrightarrow E^u(y;\ep_0)$
is well-defined and $C^2$ (with uniformly bounded derivatives) for $y\in W^u_{\ep_2}(x)$.
Moreover, there exists a constant $D_1 > 0$ such that for any such $x,y$ and $u, v\in E^u(x; \ep_1)$ we have
$\|dG_x^y(v)\| \leq D_1$, $\|G_x^y(u) - G_x^y(v)\| \leq D_1\| u-v\|$, $\|dG_x^y(u) - dG_x^y(v)\| 
\leq D_1\| u-v\|$, and  
$\|G_x^y(u) - G_x^y(v) - dG_x^y(v)\cdot (u - v)\|\leq D_1\, \|u - v \|^2$
for $\|u\|, \|v\| < \ep_1$. 

Given $x$, $y$ as above, set $x_p = f^{-p}(x)$, $y_p = f^{-p}(y)$, and notice that
\be
G_x^y \circ \hf^p_{x_p}(w) = \hf^p_{y_p} \circ G_{x_p}^{y_p}(w)
\ee
for any $w\in E^u (x_p, \ep_1)$ with $\|\hf^p_{x_p}(w)\|\leq \ep_1$. 

Let $\xi \in E^u(x;\ep_2/2)$. Setting $y = \exp^u_x(\xi)$ and $\eta = (\exp^u_y)^{-1}(x)$, we have 
$G_x^y(\xi) =  0$ and $G_x^y(0) = \eta$.
We will show that
\be
F_x(u) = d G_y^x(\eta) \circ \Lye \left( F_y\circ G_x^y (u) -  F_y (\eta)\right) \;,
\ee
for any $u \in E^u(x;\ep_2/2)$. This would imply that there exists 
$$dF_x(\xi) = d G_y^x(\eta) \circ \Lye\circ dF_y(0) \circ dG_x^y (\xi) = d G_y^x(\eta) \circ \Lye 
\circ dG_x^y (\xi)\;,$$
so $F_x$ is $C^1$ on $E^u(x;\ep_2/2)$.

To prove (3.14), consider any $u \in E^u(x; \ep_2)$, and set $v = G_x^y(u) \in E^u(y)$,
$u_p = d\hf_{x_p}^p(0) \cdot \hf^{-p}_x(u)$, $v_p = d\hf_{y_p}^p(0) \cdot \hf^{-p}_y(v)$,
$\xi_p = d\hf_{x_p}^p(0) \cdot \hf^{-p}_x(\xi)$, $\eta_p = d\hf_{y_p}^p(0) \cdot \hf^{-p}_y(\eta)$, 
$\tu_p =  \hf^{-p}_x(u)$, $\tv_p =  \hf^{-p}_y(v)$,
$\txi_p = \hf^{-p}_x(\xi)$, $\teta_p =  \hf^{-p}_y(\eta)$. Then  by (3.13), $G_{x_p}^{y_p}(\tu_p) = \tv_p$,
$G_{x_p}^{y_p}(\txi_p) = 0$ and  $G_{x_p}^{y_p}(0) = \teta_p$. Thus,
$$\tu_p =  G_{y_p}^{x_p} (\tv_p) - G_{y_p}^{x_p} (\teta_p) = dG_{y_p}^{x_p} (\teta_p) \cdot (\tv_p - \teta_p) + w'_p\;,$$
for some $w'_p$ with $\|w'_p\| \leq D_1 \, \|\tv_p - \teta_p\|^2$, and therefore, using (3.13),
\begin{eqnarray*}
u_p 
& = & d\hf^p_{x_p}(0) \cdot \tu_p = d\hf^p_{x_p}(0) \circ dG_{y_p}^{x_p} (\teta_p) \cdot (\tv_p - \teta_p) + 
d\hf^p_{x_p}(0)\cdot w'_p \\
& = &  dG_{y}^{x} (\eta) \circ d\hf^p_{y_p}(\teta_p)  \cdot (\tv_p - \teta_p) +  d\hf^p_{x_p}(0)\cdot w'_p
\end{eqnarray*}
Next, 
$d\hf^p_{y_p}(\teta_p)  \cdot (\tv_p - \teta_p) = \Lyep \circ \Fyp (v) - \Lyep \circ \Fyp (\eta) 
\to \Lye \circ F_y(v) - \Lye\circ F_y(\eta)$
as $p \to \infty$, and 
$$\|d\hf^p_{x_p}(0)\cdot w'_p\| \leq \lambda_{x_{1}}  \ldots \lambda_{x_{p}} \, \|w'_p\|
\leq D_1^2\, \frac{\lambda_{x_{1}}  \ldots \lambda_{x_{p}}}{\mu^2_{x_{1}}  \ldots \mu^2_{x_{p}}}\, 
\|u\|\, \|\eta\|
\leq D_1^2\, \gamma^p\, \|u\| \,\|\eta\| \;,$$
so $\lim_{p\to \infty} \|d\hf^p_{x_p}(0)\cdot w'_p\|  = 0$.  Thus,
$F_x(u) = \lim_{p\to \infty} u_p = dG_{y}^{x} (\eta)\circ \Lye \cdot (F_y(v) - F_y(\eta)) \;, $
which proves (3.14). Hence $F_x$ is $C^1$ on $E^u(x;\ep_2/2)$. Replacing $\ep_2$ by a smaller number we
have that $F_x$ is $C^1$ on $E^u(x;\ep_2)$.

\ms

(c) This follows easily from the definition of $F_x$.
\endofproof


\section{Ball size comparison in Bowen's metric}
\setcounter{equation}{0}

Let $M$ be a $C^2$ Riemann manifold and let  $\mt$ be a basic set for a $C^2$ flow 
$\phi_t : M \longrightarrow M$ satisfying the pinching condition (P), and let
$n$ be the {\it dimension of the (local) unstable manifolds} $W^u_\ep(x)$, $x\in \mt$.

To prove Proposition 1.2, we will first  establish the following local version.

\bs

\noindent
{\sc Lemma 4.1.} {\it  There exists a constant $\ep_3 > 0$  with the following properties:}

\ms

(a) {\it For any  $x\in \mt$ and any $0 < \delta \leq \ep \leq \ep_3$ there exist  a  constant $C =  C (x,\delta , \ep) > 0$ 
and an open neighbourhood $V_0 = V_0(x,\delta)$ of $x$ in $W^s_{\ep_0}(x) \cap \mt$ such that
$\diam \left( \mt \cap B^u_{T}(\phi_{-T}( y),\ep)  \right) \leq  C  \, \diam \left( \mt \cap B^u_T (\phi_{-T} ( y) , \delta) \right) \;$
for any $y \in V_0$ and any $T \geq  0$.}

\ms

(b) {\it For any  $x\in \mt$ and any $0 < \ep \leq \ep_3$  there exists an open neighbourhood 
$V_0 = V_0(x,\ep)$ of $x$ in 
$W^s_{\ep_0}(x) \cap \mt$  with the following property: for any $\rho \in (0,1)$ there exists 
$\delta  \in (0,\ep]$
such that for  any $y \in V_0$ and any  $T \geq 0$ we have}
$\diam \left( \mt \cap B^u_T(\phi_{-T}(y),\delta)  \right) \leq  
\rho \, \diam \left( \mt \cap B^u_T(\phi_{-T} (y) , \ep) \right) \;.$

\ms

Fix $t_0 > 0$ with (3.2).  The compactness of $\mt$ and the smoothness of the flow $\phi_t$ imply 
the existence of a constant $C_2 \geq 1$ such that
\be
d(\phi_t(y), \phi_t(z)) \leq C_2 \, d(y,z) \quad , \quad y,z\in \mt\;,  |t|\leq t_0\:.
\ee
For a non-empty set $X\subset E^u (x)$ and  $r > 0$ set
$$\ell(X) = \sup \{ \|u\| : u \in X\}\quad , \quad X(r) = \{ u\in X : \| u \| \leq r\}\;.$$

Before proving Lemma 4.1 we will first use it to derive Proposition 1.2.  

\bs

\noindent
{\sc Lemma 4.2.} {\it  There exists a constant  $\hep_0 \in (0,\ep_3]$  with the following properties:}

\ms

(a) {\it For any $x \in \mt$  and any  $0 < \delta \leq \ep \leq \hep_0$ there exist 
a constant $R_x =  R (x,\delta , \ep) > 0$ and an open neighbourhood $\oo_x$ of $x$ in $\mt$ 
such that 
$\ell \left( \mt \cap B^u_T(z,\ep)  \right) \leq R_x \, \ell \left( \mt \cap B^u_T (z, \delta) \right) \;$
for any $z \in \mt$ and $T > 0$ with $\phi_T(z) \in \oo_x$.}

\ms

(b) {\it For any $x \in \mt$,  any  $0 < \ep \leq \hep_0$ and any $\rho \in (0,1)$ there exist 
$\delta \in (0, \ep)$ and an open neighbourhood $\oo_x$ of $x$ in $\mt$  such that 
$\ell \left( \mt \cap B^u_T(z,\delta)  \right) \leq \rho \, \ell \left( \mt \cap B^u_T (z, \ep) \right) \;$
for any $z \in \mt$ and  $T > 0$ with $\phi_T(z) \in \oo_x$.} 

\bs

\noindent
{\it Proof of Lemma 4.2.}   Set $\hep_0 = \frac{\ep_3}{2C^2_2 }$.

(a) Assume that $0 < \delta, \ep \leq \hep_0$.
Let $x\in \mt$ and let $V_0 = V_0 (x,\delta/(2C_2))$ be the neighbourhood of $x$ in $\mt\cap W^s_{\ep_0}(x)$
defined as in Lemma 4.1(a) replacing $\delta$ by $\delta/(2 C_2)$. Let  $C = C (x, \delta/(2 C_2), 2C_2^2 \ep )> 0$ 
be the constant from Lemma 4.1(a) with $\delta$ and $\ep$ replaced by $\delta/(2 C_2)$ and $2C_2^2 \ep$, respectively.

Since the local product $[\cdot , \cdot ]$ and the temporal distance function $\Delta$ (see section 2) are
continuous on $\mt$, there exists an open neighbourhood $\oo_{x} = \oo_{x}(\delta, \ep)$ of $x$ in $\mt$ such that
\be
[x,z]\in V_0 \quad, \quad |\Delta (x,z)| \leq t_0
\quad, \quad d(\phi_{-\Delta(x,z)} z, [x,z]) \leq \frac{\min\{\delta, \ep\}}{2 C_2} \quad, \quad z\in \oo_x\;.
\ee

We will now check that $\oo_{x}$ and $R_{x} = C_2^2 C$  have the required properties.
Let $z \in \oo_{x}$ and $T > 0$. Set $y = [x,z]$, $t = -\Delta(x,z)$ and $\zeta = \phi_t(z)$. 
Then $y \in \mt \cap W^s_{\ep_0}(x)$ and $\phi_{-t}(y) \in W^u_{\ep_0}(z)$, so $\zeta \in W^u_{\ep_0}(y)$.
Moreover, it follows from (4.2) that $y \in V_0$, $|t| \leq t_0$ and 
\be
d(\zeta,y) \leq \frac{\min\{ \delta, \ep\}}{2 C_2}\;.
\ee
We will also need the points $z' = \phi_{-T}(z)$, $y' = \phi_{-T}(y)$ and $\zeta' = \phi_{-T} (\zeta)$.
Clearly $\zeta' = \phi_t(z') \in W^u_{\ep_0}(y')$. 

We claim that
\be
\mt \cap B^u_T(z', \ep) \subset \phi_{-t} (\mt \cap B^u_{T}(y', 2C_2^2 \ep)\;.
\ee
Indeed, let $\xi' \in \mt \cap B^u_T(z', \ep)$. Then $d(\phi_T(\xi'), \phi_T(z')) \leq \ep$, so
for $\eta' = \phi_t(\xi') \in \mt \cap W^u_{\ep_0}(\zeta')$, (4.1) implies 
$d(\phi_T(\eta'), \phi_T(\zeta')) \leq C_2\ep$, i.e. $d(\phi_T(\eta'), \zeta) \leq C_2\ep$. This and (4.3) yield
$d(\phi_T(\eta'), \phi_T(y')) = d(\phi_T(\eta'), y)  \leq 2C_2 \ep$. Since $\phi_{t_0}$ is expanding 
on local unstable manifolds of size $\ep_0$ by the choice of $t_0$, combining the latter with (4.1) gives
$d(\phi_s(\eta'), \phi_s(y')) \leq 2C_2^2 \ep$ for all $s\in [0,T]$. Thus, $\eta' = \phi_t(\xi') \in B^u_T(y', 2C_2^2\ep)$.
This proves (4.4).

Next, we will show that
\be
\phi_{-t}(\mt \cap B^u_T(y', \delta/(2C_2) ))  \subset \mt \cap B^u_{T}(z', \delta) \;.
\ee
To prove this, consider any $\eta' \in \mt \cap B^u_T(y', \delta(2C_2) )$ and let
$\xi' = \phi_{-t}(\eta')$. Then 
$$d(\phi_T(\eta'), y) = d(\phi_T(\eta'), \phi_T(y')) \leq \frac{\delta}{2C_2}\;,$$
and (4.3) implies 
$$ d(\phi_T(\eta'), \phi_T(\zeta')) = d(\phi_T(\eta'),  \zeta)  \leq \frac{\delta}{C_2}\;.$$
Since $\eta' = \phi_t(\xi')$ and $\zeta' = \phi_t(z')$, combining the latter with (4.1) gives
$d(\phi_s(\xi'), \phi_s(z')) \leq \delta$ for all $s\in [0,T]$. Thus, $\xi' \in B^u_T(z', \delta)$, which
proves (4.5).

Finally, using (4.4), (4.5) and (4.1) and with $C > 0$ defined above, one gets
\begin{eqnarray*}
\ell(\mt \cap B^u_{T}(z', \ep) )
 \leq  C_2\, \ell(\mt \cap B^u_{T}(y', 2C_2^2\, \ep) )
\leq C_2 \, C \, \ell(\mt \cap B^u_{T}(y', \delta/(2C_2)) )
 \leq  C_2^2\, C\, \ell(\mt \cap B^u_{T}(z', \delta) )\;.
\end{eqnarray*}
This proves part (a).

The proof of part (b) is very similar to the above and we omit it. \endofproof

\bs

\noindent
{\it Proof of Proposition 1.2.} Choose $\hep_0 > 0$ as in Lemma 4.2.

(a) Let $0 < \delta \leq \ep \leq \hep_0$.
It follows from Lemma 4.2(a) that for any $x \in \mt$  there exist 
a constant $R_x =  R (x,\delta , \ep) > 0$ and an open neighbourhood $\oo_x$ of $x$ in $\mt$ such that 
$\ell(\mt \cap B^u_{T}(z, \ep) ) \leq  R_x\, \ell(\mt \cap B^u_{T}(z, \delta) )$
for any $z \in \mt$ and $T > 0$ with $\phi_T(z) \in \oo_x$. Since $\mt$ is compact, there exist finitely many
neighbourhoods $\oo_{x_1}, \ldots, \oo_{x_m}$ covering $\mt$. Then
$R = 2 \max_{1\leq j\leq m} R_{x_j} > 0$ satisfies (4.1) for any $z \in \mt$ and any $T > 0$.

The proof of part (b) in the definition of regular distortion along unstable manifolds is similar 
and we omit it. 
 \endofproof

\bs

The rest of this section is devoted to the proof of Lemma 4.1. 

For $y\in \mt$, $\ep \in (0,\ep_0]$ and an integer $p \geq 0$, set
$$\hB^u_p(y , \ep) = \{ v\in E^u(y;\ep_0) :  \|\hf^p_y(v)\| \leq \ep\}\;.$$
Since  the maps $\hf_z$ are expanding distances on $E^u(z;\ep_1)$, clearly
$v \in \hB^u_p(y,\ep)$ is equivalent to $\|\hf^j_y(v)\| \leq \ep$ for all $j = 0,1,\ldots, p$.

Choose the constants $0 < \ep_2 < \ep_1 \leq \ep_0$ and $C_1 > 1$ as in  section 3, assuming that
$6 \ep_2\, C_1 <1$. In what follows we will use the notation from section 3.

Assuming $\ep_1 \in (0,\ep_0]$  is sufficiently small, for any $x\in \mt$ and 
$y \in \mt\cap W^s_{\ep_1}(x)$  the local {\it holonomy map}
$\hhs_{x,y} : W^u_{\ep_1}(x) \longrightarrow W^u_{\ep_0}(y)$ along stable laminations is 
well-defined and uniformly H\"older continuous (see section 2). 
Further, we will assume that the constant  $\ep_2 $ from Lemma 3.3 is chosen so small that for any 
$x\in \mt$ and any $y \in \mt \cap W^s_{\ep_1}(x)$ the {\it pseudo-holonomy map}
$$\hhhs_{x,y} = (\exp^u_y)^{-1}\circ \hhs_{x,y} \circ \exp^u_x :  
E^u (x;\ep_2) \longrightarrow E^u (y;\ep_1)$$
is well-defined and uniformly H\"older.  
Notice that    
$\hhhs_{x,y}( \hLa_x(\ep_2)) \subset \hLa_y(\ep_1)$ for any  
$x\in \mt$ and any $y \in \mt\cap W^s_{\ep_1}(x)$.

Thus, instead of dealing with sets of the form $\mt \cap B^u_T(\phi_{-T}(y), \ep)$ in Lemma 4.1, 
it is enough to prove the analogous statements for sets of the form  
$\mt \cap B^u_{pt_0}(f^{-p}(y), \ep)$, $p\geq 1$, which in turn
combined with the local uniform Lipschitzness of the maps $\exp^u_y$ leads to analogous statements
for sets of the form $\hLa_{f^{-p}(y)}(\ep) \cap \hB^u_{p}( f^{-p}(y),\ep)$.

Recall the maps $F_x$ from Theorem 3.1. For the proof of Lemma 4.1 we will also need the sets
$$\tB^u_p (z, \ep) = \{ u \in E^u (z) : \|d \hf^p_z(0) \cdot u\| \leq\ep\}\quad , \quad 
\tLa_x = F_x(\hLa_x(\ep_2)) \subset E^u (x; 2\ep_2)\;,$$
and the maps $\thhs_{x,y} = F_y \circ \hhhs_{x,y} \circ (F_x)^{-1} :  
E^u (x;\ep'_2) \longrightarrow E^u (y;\ep_1)\;.$
Clearly we can take $\ep'_2 \in (0,\ep_2]$ independent of $x$ and $y$ so that the above 
is well-defined and uniformly 
H\"older for any $x\in \mt$ and any $y \in \mt \cap W^s_{\ep_1} (x)$. 
Moreover we have $\thhs_{x,y}(\tLa_x(\ep'_2)) \subset \tLa_y(\ep_2)\;.$
Another property of the sets $\tLa_x$ is contained in the following immediate consequence of Theorem 3.1.

\bs

\noindent
{\sc Corollary 4.3.} {\it For any $x \in \mt$, any $\ep \in (0,\ep_2/2]$  and any $p\geq 1$, 
setting $x_p = f^{-p}(x)$ we have 
$d\hf^p_{x_p}(0) (\tLa_{x_p}(\ep)\cap \tB^u_p(x_p, \ep)) \subset \tLa_x(\ep)\;.$
More generally, $d\phi_t(x)\cdot \tLa_x(\ep) \subset \tLa_{\phi_t(x)}(\ep)$ for all 
$x\in \mt$$, t\leq 0$ and $\ep\in (0,\ep_2]$.}

\bs

The following lemma is rather important for the proof of the central Lemma 4.1. It describes some sort of a 
`tangent bundle' $E^u_\mt(x)$ ($x \in \mt$) to the set $\mt$ which is $d\phi_t$-invariant and
has  some continuity properties, as well. 

Given $x\in \mt$, let $m_x \geq 1$ be the minimal integer such that there exists $\ep(x) \leq \ep_2$
with $\dim (\span (\tLa_x(\delta))) = m_x$ for all $0 < \delta \leq \ep(x)$. Then the linear subspace
$E^u_\mt(x) = \span (\tLa_x(\delta))$ is the same for all $\delta \in (0,\ep(x)]$. Corollary 4.3 shows that
$m_{\phi_t(x)} = m_x$ and $d\phi_t(x)(E^u_\mt (x)) = E^u_\mt (\phi_t(x))$ for all $x\in \mt$ and $t\in \R$.

\bs

\noindent
{\sc Lemma 4.4.} {\it There exists an integer $m$ such that $m_x = m$ for any $x \in \mt$. Moreover,
for any $x\in \mt$ we have $E^u_\mt (x) = \span (\tLa_x(\ep'_2))$, where $\ep'_2 > 0$ is as above, 
and there exists 
$\ep = \ep(x)\in (0,\ep_0)$ such that $E^u_\mt(y)$ depends continuously on $y \in W^s_\ep(x) \cap \mt$.}

\bs

\noindent
{\it Proof of Lemma} 4.4. 
Set $m = \min_{x\in \mt} m_x$ and let $y\in \mt$ be such that $m_{y} = m$. 
Let $\ep = \ep(y) \in (0,\ep_2]$ be small enough so that
$E^u_\mt (y) = \span(\tLa_{y}(\delta))$ for all $\delta \in (0,\ep]$. Let $v_1, \ldots, v_{m} \in \tLa_y(\ep)$
be a linear basis in $E^u_\mt(y)$. Assume that $0 < \delta < \min\{ \ep_2/6, \ep/6\}$, where 
$\ep_2\in (0,\ep_1/(2C_1)]$
is as in the proof of Theorem 3.1. We will now use  the map
$G_x^y = (\exp^u_y)^{-1}\circ \exp^u_x : E^u(x;\ep_1) \longrightarrow E^u(y;\ep_0)$ 
and formula (3.10) from that proof.

Let $x\in \mt\cap W^u_\delta(y)$. We will show that $m_x \leq m$, so we must have $m_x = m$. 
Setting $\eta = (\exp^u_{y})^{-1}(x)$, (3.14) holds for any $u \in E^u(x;\ep_2)$. Clearly 
$\eta \in \hLa_y(\delta)$, so
$F_y(\eta)\in \tLa_y(2\delta)$ (since by Theorem 3.1 (b), $\|F_y(\eta)\| \leq 2\|\eta\|$). 
Thus, $F_y(\eta)$ is a 
linear combination of the vectors $v_1, \ldots, v_m$. Given $u' \in \tLa_x(\delta)$, 
we have $u' = F_x(u)$
for some $u \in \hLa_x(2\delta)$. Then $G_x^y(u)\in \hLa_y(3\delta)$, so 
$F_y(G_x^y(u)) \in \tLa_y(6\delta) \subset \tLa_y(\ep)$, and therefore 
$F_y(G_x^y(u))$ is a linear combination of the vectors $v_1 ,\ldots, v_m$. Thus, 
$u' = F_x(u) = d G_y^x(\eta) \circ \Lye \left( F_y\circ G_x^y (u) -  F_y (\eta)\right) $
is a linear combination of the
vectors $w_j = dG_y^x(\eta) \circ \Lye \cdot v_j $ ($ j = 1, \ldots, m$), so 
$\dim(\span(\tLa_x(\delta))) \leq m$.
Hence $m_x \leq m$, and therefore $m_x = m$. Moreover, it follows from this argument that
$\span(\tLa_x(\delta)) = E^u_\mt(x)$ for any $x \in  \mt\cap W^u_\delta(y)$ and any 
$0 < \delta < \min\{ \ep_2/6, \ep/6\}$.

We now claim that $m_z = m$ for any $z\in \mt$. Assume that $m_z > m$ for some $z \in \mt$, and take 
$\ep' \in (0,\ep_2]$ so small that $\span (\tLa_z(\ep'')) = E^u_\mt(z)$ for all $\ep''\in (0,\ep']$.  
Let $u_1, \ldots, u_{m_z} \in \tLa_z(\ep')$ be a linear basis in $E^u_\mt(z)$. 
Take $\mu \in (0, \ep'_2]$ so small that
for any $z'\in W^s_\mu(z) \cap \mt$ the vectors $u_j(z') = \thhs_{z,z'}(u_j)$ 
($j = 1, \ldots, m_z$) are linearly independent 
in $E^u(z')$. By Corollary 4.3 these vectors belong to $\tLa_{z'}(\ep_2)$. 
Let $y\in \mt$ and $\ep = \ep(y) > 0$
be as above and let $0 < \delta < \min\{ \ep_2/6, \ep/6\}$. It is well-known 
(see e.g. \cite{kn:KH}) that if $T > 0$ is
sufficiently large, then $\phi_t(W^u_\delta(y)) \cap W^s_\mu(z) \neq \e$ for any $t \geq T$. 
Take $T > 0$ with this property
so that $e^{\alpha T}/C > \ep_2$. Then for some $x \in W^u_\delta(y) \cap \mt$ and some $t \geq T$ we have
$z' = \phi_t(x)\in W^s_\mu(z) \cap \mt$.  The choice of $T$ and  $t \geq T$ imply 
$d\phi_t(x) (E^u(x;\delta)) \supset E^u(z';\ep_2)$,
so by Corollary 4.3, $d\phi_t(x) (\tLa_x(\delta)) \supset \tLa_{z'}(\ep_2)$.
Since $\dim(\span(\tLa_x(\delta))) = m$, we now get $\dim(\span(\tLa_{z'}(\ep_2))) = m$, a contradiction
with the linear independence of the vectors $u_j(z') = \thh_z^{z'}(u_j)$ ($j = 1, \ldots, m_z$) and $m_z > m$.
Thus, $m_z = m $ for all $z \in \mt$. 

Using the above notation (with $m_z = m$), by the previous argument, for any 
$z'\in W^s_\mu(z) \cap \mt$ the vectors $u_j(z') = \thhs_{z, z'}(u_j)$ ($j = 1, \ldots, m$) provide 
a basis for 
$E^u_\mt(z')$, so the latter depends continuously on $z'$. Moreover, repeating the above argument we
can see that $\dim(\span (\tLa_z(\ep'_2))$ cannot exceed $m$, so we must have
$\span (\tLa_z(\ep'_2)) = E^u_\mt(z)$. The same argument can be applied to any 
$z'' \in W^s_\mu(z)\cap \mt$.
\endofproof

\bs

\noindent
{\it Proof of Lemma 4.1.}  Notice that if $x \in \mt$ and $z = x_p = f^{-p}(x)$ for some $p \geq 0$, 
then for any $\ep \in (0,\ep_2]$ we have
\be
\tLa_z(\ep/2) \cap \tB^u_p( z , \ep/2) \subset F_z(\hLa_z(\ep) \cap \hB^u_p( z , \ep))
\subset \tLa_z(2\ep) \cap \tB^u_p( z , 2\ep)\;.
\ee

(a) Let $x\in \mt$. Given $y \in \mt$ and an integer $p \geq 0$, set $y_p = f^{-p}(y) \in \mt$. 
According to (4.6), it is enough to prove the following

\bs

\noindent
{\sc Sublemma 4.5.} {\it For any  $0 < \delta \leq \ep \leq \ep'_2/2$ there exist  a  constant 
$D =  D (x,\delta , \ep) > 0$ 
and an open neighbourhood $V_0$ of $x$ in $W^s_{\ep_0}(x) \cap \mt$ such that
$\ell \left( \tLa_{y_p}(\ep) \cap \tB^u_p( y_p,\ep)  \right) \leq 
D  \, \ell \left( \tLa_{y_p}(\delta) \cap \tB^u_p ( y_p , \delta) \right) \;$
for any $y\in V_0$ and any integer $p \geq 0$.}

\ms

\noindent
{\it Proof of Sublemma 4.5.} Choose $\ep(x) \in (0,\ep'_2]$ so that $E^u_\mt(y)$
depends continuously on $y \in W^s_{\ep(x)} (x)\cap \mt$. For any $y \in W^s_{\ep(x)}(x) \cap \mt$
choose and fix an orthonormal basis $e_1(y) ,e_2(y) , \ldots,e_m(y)$ in $E^u_\mt (y)$
which depends continuously on $y$. 

Let $0 < \delta \leq  \ep \leq \ep'_2/2$. 
By the definition of $E^u_\mt(x)$ and Lemma 4.4, there exist 
$u_1, u_2, \ldots, u_m \in \tLa_x(\delta/2)$ which are 
linearly independent. Set
$\Delta = \Delta(x,\delta) = \Vol_m[u_1,u_2, \ldots,u_m] > 0\;,$
where $[u_1,u_2, \ldots,u_m] $ denotes the parallelepiped in $E^u_\mt(x)$ determined by the vectors
$u_1, \ldots,u_m$, and $u_j(y) = \thhs_{x,y}(u_j)$ for any $j = 1, \ldots,m$.
Choose an  open neighbourhood $V_0$ of $x$ in $W^s_{\ep(x)}(x) \cap \mt$ such that
\be
\Vol_m[u_1(y),u_2(y), \ldots,u_m(y)] \geq \frac{\Delta}{2}  \quad , \quad y \in V_0\;,
\ee
and
\be
\frac{\|u_j\|}{2} \leq \| u_j(y)\| \leq 2\|u_j\| \quad, \quad y\in V_0\:, \: 1\leq j\leq m\;. 
\ee
Then $u_j(y) \in \hLa_y(\delta)$ for all $j = 1, \ldots,m$.
Let $L_y = L(x,y,\delta) :  E^u_\mt(y) \longrightarrow E^u_\mt(y)$ be the linear operator 
such that  $L_y u_j(y) = e_j(y)$ 
for all $j = 1, \ldots,m$. It follows from (4.7) and (4.8) that there exists a constant
$b = b(x,\delta) > 0$ (determined by $\Delta$ and $\| u_1\|, \ldots, \|u_m\|$) such that
$\|L_y\| \leq b$ for all $y\in V_0$. 
 
Fix for a moment $y \in V_0$. Consider an arbitrary integer $p \geq 1$ and set $z = f^{-p}(y) \in \mt$. Given
$v \in \tLa_z(\ep) \cap \tB^u_p( z , \ep)$, we have $\|d\hf_z^p(0)\cdot v\| \leq \ep \leq \ep_2/2$, so
by Corollary 4.3 and Lemma 4.4,  we have $u = d\hf_z^p(0)\cdot v \in \tLa_y(\ep) \subset E^u_\mt(y)$. 
Consequently, $u = \sum_{s=1}^m c_s\, u_s(y)$ for some real numbers $c_s$, so
$L_y u = \sum_{s=1}^m c_s L_y (u_s(y))  = \sum_{s=1}^m c_s e_s(y) \;.$
Thus,  $\sqrt{\sum_{s=1}^m c_s^2} = \|L_y u\| \leq \|L_y \|\, \|u\| \leq \ep\, b \;,$ and so  $|c_s| \leq \ep\, b$ for all $ s = 1, \ldots, m$. 
Since $v_j = d\hf_y^{-p}(0) \cdot u_j(y) \in \tLa_{z}(\delta) \cap \tB^u_p(z,\delta)$ for all $j = 1,  \ldots,m$, it now follows that
\begin{eqnarray*}
\|v\|  =  \|d \hf^{-p}_y (0) \cdot u\| = \left\| \sum_{s=1}^m c_s\,  d\hf^{-p}_y(0) \cdot u_s(y)\right\|
 \leq  m \, \ep \, b \, \max_{1\leq s\leq m} \|v_s\| \leq  m \, \ep \, b \, \ell(\tLa_z(\delta) \cap \tB^u_p(z,\delta))\;.
 \end{eqnarray*}
Hence 
$\ell \left( \tLa_{z}(\ep) \cap \tB^u_p( z ,\ep)  \right)  
\leq D \, \ell \left( \tLa_{z}(\delta) \cap \tB^u_p ( z , \delta) \right)$,
where $D = D(x,\delta, \ep) = m\, \ep \, b$.
This concludes the proof of the Sublemma and thus the proof of part (a) in Lemma 4.1.

\bs

\noindent
{\it Proof of Lemma 4.1(b)}. 
We will essentially repeat the argument in the proof of Sublemma 4.5.
As before, it is enough to prove the analogous statement for sets of the form
$\tLa_{f^{-p}(x)}(\ep) \cap \tB^u_p( f^{-p}(x),\ep)$.

Choose $\ep(x)$ as before and for each $y \in W^u_{\ep(x)}(x) \cap \mt$ an orthonormal basis 
$\{ e_j(y)\}$ in $E^u_\mt (y)$ 
depending continuously on $y$.
Let $0 < \ep  \leq \ep'_2/2$ and $\rho \in (0,1)$. As in the proof of the Sublemma,
there exists a basis $u_1, u_2, \ldots, u_m$ in $E^u_\mt(x)$ with $u_j \in \tLa_x(\ep/2)$ for 
all $j = 1, 2, \ldots,m$.
Let $L_y = L(x,\ep) :  E^u_\mt(y) \longrightarrow E^u_\mt(y)$ be the linear operator such that  
$L_y u_j(y) = e_j(y)$ 
for all $j = 1, \ldots,m$, where $u_j(y)$ are defined as in the proof of the Sublemma. 
Choose $V_0$  and $b = b(x,\ep) > 0$ as
before, replacing $\delta$ by $\ep$. Let  $0 < \delta \leq  \min\left\{ \ep , \frac{\rho}{m\, b}\right\}.$ 
Then given $y \in V_0$ and an  integer $p \geq 1$, set $z = f^{-p}(y) \in \mt$. As in the proof of 
the Sublemma, for
$v \in \tLa_z(\delta) \cap \tB^u_p( z , \delta)$ one obtains 
$\|v\| \leq  m \, \delta \, b\, \ell(\tLa_z(\ep) \cap \tB^u_p(z,\ep))$.
Thus,
$ \di \ell(\tLa_z(\delta) \cap \tB^u_p(z,\delta)) \leq \rho \, \ell(\tLa_z(\ep) \cap \tB^u_p(z,\ep))$.
\endofproof

\def\chBo{\check{B}^{u,1}}
\def\tBo{\tB^{u,1}}
\def\hBo{\hB^{u,1}}
\def\hpi{\hat{\pi}}

\section{Proof of Theorem 1.1}
\setcounter{equation}{0}

Let again $M$ be a $C^2$ complete Riemann manifold, $\phi_t$  be a  $C^2$  flow on $M$, and let 
$\mt$ a basic set for $\phi_t$. Throughout we assume that $\phi_t$ and $\mt$ satisfy the conditions
(LUPC) and (I) from the Introduction.

The proof of Theorem 1.1 is a generalization of what we did 
under the pinching condition (P) in sections 3 and 4. 
As before, given $x\in \mt$ and $\ep > 0$, we have to deal with 
diameters of sets of the form $f^{-p}(\mt \cap B(x,\ep))$ ($p\geq 1$), where 
$f = \phi_{t_0}$ for some sufficiently 
large $t_0 > 0$.  Obviously, going backwards along the flow, the greatest expansion
occurs in the direction of vectors in $E^u_1$, so the diameter of $f^{-p}(\mt \cap B(x,\ep))$ would be
comparable with that of its `projection' $\pi^{u,1}(f^{-p}(\mt \cap B(x,\ep)))$ to the 
corresponding leaf of $W^{u,1}$ (see Lemma 5.1 below). 
The behaviour of $d\phi_t$ on $E^u_1$ (and that of $\phi_t$ on $W^{u,1}$) is very similar to 
what we had in sections 3 and 4, and we use the arguments from there to compare diameters of sets of the form
$\pi^{u,1}(f^{-p}(\mt \cap B(x,\ep)))$.

We now proceed with the proof of Theorem 1.1. 

Notice that, since the splitting $E^u_1(x) \oplus E^u_2(x)$ depends continuously on $x$ and $\mt$ is compact,
the angle between $E^u_1(x)$ and $E^u_2(x)$ is uniformly bounded below by a positive constant. So, the local
submanifolds $W^{u,1}_{\ep_0}(x)$ and $W^{u,2}_{\ep_0}(x)$ of $W^u_{\ep_0}(x)$ are (uniformly) 
transversal. Moreover, taking $\ep_1 \in (0,\ep_0]$ and $\ep_2 \in (0,\ep_1]$ sufficiently small, 
for any $y, z \in W^u_{\ep_1}(x)\cap \mt$ with
$d(y,z) \leq \ep_2$, the submanifolds $W^{u,1}_{\ep_1}(y)$ and $W^{u,2}_{\ep_1}(z)$ of 
$W^u_{\ep_0}(x)$ are transversal 
and of complementary dimension, and intersect at a single point 
$[y,z]_x^u = W^{u,1}_{\ep_1}(y) \cap W^{u,2}_{\ep_1}(z)\;.$
It follows immediately that the so defined local product has the usual invariance, namely
$\phi_t([y,z]_x^u) = [\phi_t(y), \phi_t(z)]_{\phi_t(x)}^u$ for $t \leq 0$. However, {\bf in general $[y,z]_x^u$ does not have to belong to} $\mt$.
We can now define the {\it projection 
$\pi^{u,1}_x : W^u_{\ep_2}(x)\cap \mt \longrightarrow W^{u,1}_{\ep_1}(x)$ along
$W^{u,2}$} by $\pi^{u,1}_x(y) = [x,y]^u_x$. 

In what follows for any $u \in E^u(x)$, $x\in \mt$, we will use the notation 
$u = u_1+u_2$, where $u_i\in E^u_i(x)$
for $ i = 1,2$. Setting $\|u\|' = \max\{ \|u_1\|, \|u_2\|\}$ defines a norm on 
$E^u(x)$ equivalent to the original norm
$\|u\|$ defined by the Riemann metric on $M$. For a non-empty subset $X$ of $E^u(x)$ let $\diam'(X)$ be the
{\it diameter} of $X$ with respect to $\|\cdot \|'$.

Fix an arbitrary and sufficiently large $t_0 > 0$ as in section 3 and  set $f = \phi_{t_0}$. 
Then $f$ is a partially  hyperbolic diffeomorphism with respect to the invariant splitting 
$E(x) = E^s(x)\oplus (E^0(x)\oplus E^u_1(x)) \oplus E^u_2(x)$ with $E^c(x) = E^0(x)\oplus E^u_1(x)$, 
and it follows from Theorem $A'$ in \cite{kn:PSW} that the local holonomy maps along the lamination 
$W^{u,2}$ through $\mt$ are $\theta$-H\"older for some sufficiently small $\theta > 0$. In particular, 
the projections $\pi^{u,1}_x$ are uniformly continuous, and it follows from this that if $\{x_m\}$ 
and $\{y_m\}$ are sequences in $\mt$ with $y_m\in W^u_{\ep_2}(x_m)$
for all $m$, $x_m\to x \in \mt$ and $y_m \to y\in W^u_{\ep_2}(x)$ as $m \to \infty$, 
then $\pi^{u,1}_{x_m}(y_m) \to \pi^{u,1}_x(y)$ as $m \to \infty$. 

Assuming that $\ep_0 > 0$ is sufficiently small, for each $x\in \mt$ there exists a $C^2$ diffeomorphism
$\Phi_x : E^u(x;\ep_0) \longrightarrow W^u_{\ep_0}(x)$ such that 
$(\Phi_x)_{| E^u_i(x;\ep_0)} : E^u_i(x;\ep_0) \longrightarrow  W^{u,i}_{\ep_0}(x)$ is the corresponding 
exponential map for $i = 1,2$. We can choose $\Phi_x$ in such a way that the diffeomorphism
$(\exp^u_x)^{-1}\circ \Phi_x$ has uniformly bounded derivatives, so in particular if $C > 0$ 
is sufficiently large, then\\
$\frac{1}{C}\,\|u-v\| \leq  d(\Phi_x(u),\Phi_x(v)) \leq C\, \|u-v\|$ for all $x\in \mt$, 
$u,v\in E^u(x;\ep_0)$. Moreover, since
the leaves of the distribution $W^{u,2}$ are $C^1$ in $W^u$,  assuming again that $\ep_0> 0$, 
$\ep_1\in (0,\ep_0]$ 
and $\ep_2\in (0,\ep_1]$ are  sufficiently small and the constant $C > 0$ is sufficiently large, 
for every $x\in \mt$ and every $y\in \mt \cap W^{u}_{\ep_2}(x)$  we have
\be
\frac{\|v_1\|}{C} \leq  d(\pi^{u,1}_x(y),x) \leq C \|v_1\| \quad , \quad v = (v_1,v_2) = (\Phi_x)^{-1}(y)\;.
\ee

Next, assuming that $\ep_1 \in (0,\ep_0]$ is sufficiently small for any $x\in \mt$ the map
$$\hf_x = (\Phi_{f(x)})^{-1} \circ f\circ \Phi_x : E^u (x ; \ep_1) \longrightarrow E^u (f(x); \ep_0)\;$$
is well-defined and therefore $C^2$.  It is important to notice that
$$\hf_x(E^u_i (x ; \ep_1)) \subset E^u_i (f(x); \ep_0) \quad , \quad i = 1,2.$$

As in section 3, for any $y \in \mt$ and any integer $k \geq 1$ we will use the notation
$$\hf_y^k = \hf_{f^{k-1}(y)} \circ \ldots \circ \hf_{f(y)} \circ \hf_y\quad,
\quad \hf_y^{-k} = (\hf_{f^{-k}(y)})^{-1} \circ \ldots \circ (\hf_{f^{-2}(y)})^{-1} 
\circ (\hf_{f^{-1}(y)})^{-1} \;,$$
at any point where these sequences of maps are well-defined. Finally, for $x\in \mt$ and 
$\ep \in (0,\ep_0]$ set
$$\hLa_x(\ep) = \{ u\in E^u(x;\ep) : \Phi_x(u) \in \mt\}\;.$$
As before we have   $\hf^{-1}_x (\hLa_x(\ep)) \subset \hLa_{f^{-1}(x)}(\ep)$.

\ms

\noindent
{\sc Lemma 5.1.} {\it For any $\ep \in (0,\ep_1]$ there exists $\omega_\ep > 0$ such that for every 
$x\in \mt$ there exists $u \in \hLa_x(\ep)$ with $\| u_1\| \geq \omega_\ep$.}

\bs

\noindent
{\it Proof.} Let $\ep\in (0,\ep_1]$. According to (5.1), it is enough to show that there exists 
$\omega'_\ep > 0$ such that
for every $x\in \mt$ there exists $y\in \mt\cap W^u_{\ep_1}(x)$ with 
$d(x, \pi^{u,1}_x(y)) \geq \omega'_\ep$.
Assuming this is not so,  for every integer $m \geq 1$ there exists $x_m \in \mt$ with
$d(x_m, \pi^{u,1}_{x_m}(y)) \leq 1/m$ for all  $y\in \mt\cap W^u_{\ep}(x_m)$.
We may assume $x_m \to x\in \mt$ as $m \to \infty$.
The condition (I) implies that there exists $y\in \mt\cap W^u_{\ep/2}(x) \setminus W^{u,2}_{\ep/2}(x)$. 
Then $z = \pi^{u,1}_x(y) \neq x$. Setting $y_m = \hhs_{x,x_m}(y)$, where $\hhs_{x,x_m}$ is the 
local stable holonomy map (see section 2), we have $y_m\in \mt \cap W^u_{\ep}(x_m)$ for all
sufficiently large $m$ and $y_m \to y$ as $m \to \infty$. Hence 
$\pi^{u,1}_{x_m}(y_m) \to \pi^{u,1}_x(y) = z$, so for all sufficiently large $m$ we have
$d(\pi^{u,1}_{x_m}(y_m), x_m) > d(z,x)/2 > 0$. This is a contradiction with
$d(\pi^{u,1}_{x_m}(y_m), x_m) \leq 1/m$ implicated by the choice of $x_m$.
\endofproof

\bs

Next,  for any $y \in \mt$, $\ep\in (0,\ep_2]$ and $p \geq 1$ set
$$\chBo_p(y,\ep) = \{ v_1 \in E^u_1(y;\ep) : \exists 
v = (v_1,v_2)\in \hLa_y(\ep) \: \mbox{\rm with }\: \|\hf^p_y(v)\|' \leq \ep\}\;,$$
and $\hB^u_p(y,\ep) = \{ v \in E^u(y;\ep) :  \|\hf^p_y(v)\|' \leq \ep \}$.
Clearly, 
\be
\diam'(\chBo_p(y,\ep)) \leq \diam' (\hB^{u}_p(y,\ep) \cap \hLa_y(\ep)) \;.
\ee

The following consequence of Lemma 5.1 is derived by using some well-known arguments (see e.g.
Appendix A.1 in \cite{kn:BR}).

\bs

\noindent
{\sc Lemma 5.2.} {\it Choosing $\ep_1 > 0$ sufficiently small,
for any $\ep \in (0,\ep_1]$ there exists an integer $p_\ep \geq 1$ such that 
\be
\diam' (\hB^{u}_p(y,\ep) \cap \hLa_y(\ep)) \leq \diam'(\chBo_p(y,\ep)) 
\ee
for every $y\in \mt$ and every integer $p \geq p_\ep$.}

\bs

\noindent
{\it Proof.} 
Fix constants $\lambda_1 < \mu_2$ such that $e^{\beta\, t_0} < \lambda_1 < \mu_2 < e^{\alpha_2\, t_0}$. 
Using Taylor's formula, there exists a  constant $D > 0$ such that
\be
\|\hf^{-1}_x(u) - \hf^{-1}_x(v) - d\hf^{-1}_x(v)\cdot (u-v)\| \leq D\, \|u-v\|^2\quad ,\quad
x\in \mt\:, \: u, v\in E^u(x;\ep_0)\;.
\ee

Take $\ep_1 \in (0,\ep_0]$ such that 
$D\ep_1 < \min\{ 1/\mu_2 - e^{-\alpha_2 t_0}\; ,  \; e^{-\beta t_0} - 1/\lambda_1\}$. 
Given $\ep \in (0,\ep_1]$, let $\omega_\ep > 0$ be as in Lemma 5.1 and let
$p_\ep \geq 1$ be the least integer so that $(\mu_2/\lambda_1)^{p_\ep} \geq \ep/\omega_\ep$.

Setting $\hf^{-1}_x(u) = ( (\hf^{-1}_x)^{(1)}(u) , (\hf^{-1}_x)^{(2)}(u) )$, for any $x\in \mt$ and 
$u \in E^u(x;\ep_0)$, (5.4) implies
$$\hf^{-1}_x(u) - \hf^{-1}_x(u_1,0) =  d\hf^{-1}_x(u_1,0)\cdot (0,u_2) + v\;,$$
where $\|v\| \leq D\, \|u_2\|^2$. Comparing the $E^u_2$-coordinates and using the choice of $\ep_1$ gives
$$\|(\hf^{-1}_x)^{(2)}(u)\| \leq \|d\hf^{-1}_x(u_1,0)\cdot (0,u_2)\| + \|v_2\|
\leq \frac{\|u_2\|}{e^{\alpha_2 t_0}} + D\|u_2\|^2 \leq \frac{\|u_2\|}{\mu_2} \;.$$
In a similar way one gets
$\di \|(\hf^{-1}_x)^{(1)}(u)\| \geq \frac{\|u_1\|}{\lambda_1}\;.$

Let $y\in \mt$ and $p\geq p_\ep$ be an integer.
Let $w = (w_1,w_2) \in \hB^u_p(y,\ep) \cap \hLa_y(\ep)$ be such that $\|w\|'$ is the maximal possible. 
Since $w_1 \in \chBo_p(y,\ep)$,  if $\|w_1\| \geq \|w_2\|$, then $\|w_1\| \geq \|w\|'$, and so (5.3) is 
trivially satisfied in this case.

Assume that $\|w_1\| < \|w_2\|$; then $\|w\|' = \|w_2\|$. Set $x = f^p(y)$ and  $\zeta = \hf_y^p(w)$; then 
$\zeta \in \hLa_x(\ep)$. By Lemma 5.1 there exists $u \in \hLa_x(\ep)$ with $\|u_1\| \geq \omega_\ep$. 
Now $\|\zeta\| \leq \ep$ implies $\|u_1\| \geq (\omega_\ep/\ep)\, \|\zeta_2\|$,
and it follows from above that 
$$\|(\hf_x^{-1})^{(1)}(u)\| \geq \frac{\|u_1\|}{\lambda_1} \geq \frac{(\omega_\ep/\ep) 
\|\zeta_2\|}{\lambda_1} \geq 
(\omega_\ep/\ep) \, \left(\frac{\mu_2}{\lambda_1}\right)\, \|(\hf_x^{-1})^{(2)}(\zeta)\|\;. $$
Using this argument by induction, for $v = \hf^{-p}_x(u) \in \hB^u_p(y,\ep) \cap \hLa_y(\ep)$ 
we get 
$$\|v_1\| \geq (\omega_\ep/\ep)\, (\mu_2/\lambda_1)^p\, \|w_2\| \geq \|w_2\|\;.$$
Thus, $\diam'(\chBo_p(y,\ep)) \geq \|v_1\| \geq \|w\|' = \diam'(\hB^u_p(y,\ep) \cap \hLa_y(\ep))$.
\endofproof

\bs

To prove Theorem 1.1, it remains to compare diameters of sets of the form $\chBo_p(y,\ep)$. 
 As in section 4, the main step
is the following lemma whose proof follows  the arguments from sections 3 and 4 with some small 
modifications. 
For completeness we sketch its proof  in the Appendix omitting most of the details.

\bs  

\noindent
{\sc Lemma 5.3.} {\it There exists a constant $\ep_3 \in (0,\ep_2]$  with the following properties:}

\ms

(a) {\it For any  $x\in \mt$ and any $0 < \delta \leq  \ep \leq \ep_3$ there exist  a  constant 
$R =  R (x,\delta , \ep) > 0$ 
and an open neighbourhood $V_0 = V_0(x,\delta)$ of $x$ in $W^s_{\ep_0}(x) \cap \mt$ such that
\be
\diam \left( \chBo_{p}(f^{-p}( y),\ep)  \right) 
\leq  R  \, \diam \left( \chBo_p (f^{-p} ( y) , \delta) \right) \;
\ee
for any $y \in V_0$ and any integer $p \geq  1$.}

\ms

(b) {\it For any  $x\in \mt$ and any $0 < \ep \leq \ep_3$  there exists an open neighbourhood 
$V_0 = V_0(x,\ep)$ of $x$ in 
$W^s_{\ep_0}(x) \cap \mt$  with the following property: for any $\rho \in (0,1)$ there exists 
$\delta  \in (0,\ep]$
such that for  any $y \in V_0$ and any  integer $p \geq 1$ we have}
$\diam \left( \chBo_p(f^{-p}(y),\delta)  \right) \leq  \rho \, 
\diam \left( \chBo_p(f^{-p} (y) , \ep) \right) \;.$

\ms

\noindent
{\it Proof of Theorem 1.1.} As in section 4, we first derive the existence of
 a constant  $\hep_0 \in (0,\ep_3]$  with the following properties:

\ms

(i) For any $x \in \mt$  and any  $0 < \delta \leq  \ep \leq \hep_0$ there exist 
a constant $R_x =  R (x,\delta , \ep) > 0$ and an open neighbourhood $\oo_x$ of $x$ in $\mt$ 
such that 
$\ell \left( \chBo_T(z,\ep)  \right) \leq R_x \, \ell \left(\chBo_T (z, \delta) \right) \;$
for any $z \in \mt$ and $T > 0$ with $\phi_T(z) \in \oo_x$. Here
$$\chBo_T(z,\ep) = \{ v_1 \in E^u_1(z;\ep) : \exists v = (v_1,v_2)\in \hLa_z(\ep) \: 
\mbox{\rm with }\: \|\Phi_{\phi_T(z)}^{-1}\circ \phi_T\circ \Phi_z(v)\|' \leq \ep\}\;.$$

\ms

(ii) For any $x \in \mt$,  any  $0 < \ep \leq \hep_0$ and any $\rho \in (0,1)$ there exist 
$\delta \in (0, \ep)$ and an open neighbourhood $\oo_x$ of $x$ in $\mt$  such that 
$\ell \left( \chBo_T(z,\delta)  \right) \leq \rho \, \ell \left( \chBo_T (z, \ep) \right) \;$
for any $z \in \mt$ and  $T > 0$ with $\phi_T(z) \in \oo_x$.

\ms

Next, fix $\hep_0 > 0$ as above.

\ms

(a) Let $0 < \delta \leq \ep \leq \hep_0$. It follows from (i) above that for any $x \in \mt$  there exist 
a constant $R_x =  R (x,\delta , \ep) > 0$ and an open neighbourhood $\oo_x$ of $x$ in $\mt$ such that 
$\ell(\chBo_{T}(z, \ep) ) \leq  R_x\, \ell(\chBo_{T}(z, \delta) )$
for any $z \in \mt$ and $T > 0$ with $\phi_T(z) \in \oo_x$. Since $\mt$ is compact, there exist finitely many
neighbourhoods $\oo_{x_1}, \ldots, \oo_{x_m}$ covering $\mt$. Then
$R = 2 \max_{1\leq j\leq m} R_{x_j} > 0$ satisfies (5.5).

The proof of part (b) in the definition of regular distortion along unstable manifolds is similar and we omit it. 
 \endofproof

\section{Appendix: Proof of Lemma 5.3}
\setcounter{equation}{0}

We will use the notation from section 5.
Clearly, what enables us to use the arguments from sections 3 and 4 is the pinching condition 
on the spectrum of $d\phi_t$ over the bundle $E^u_1(x)$, and also the invariance of 
$E^u_1(y)$, $[y,z]^u_y$ and
$\hLa_y(\ep)$ under $\hf_y^{-1}$. 

Set
$$\Lao_y(\ep) = \{ \pi^{u,1}_y(z) : z\in \mt \cap W^u_\ep(y)\} \subset W^{u,1}_\ep(y) \quad ,\quad
\hLao_y = (\Phi_y)^{-1}(\Lao_y(\ep_2))\;.$$
It is important properties to notice that $f^{-1}(\Lao_y(\ep)) \subset \Lao_{f^{-1}(y)}(\ep)$ and
\be
\hf_y^{-1}(\hLao_y(\ep)) \subset \hLao_{f^{-1}(y)}(\ep)\;.
\ee

Notice that a set of the form $\chBo_p(y,\ep)$ is not necessarily a subset of $\hLa_y(\ep)$, however
it is contained in $\Phi_y^{-1}(\Lao_y(\ep))$.
For $0 < \ep \leq \ep_2$, $y\in \mt$ and $p \geq 0$, the set
$$\hBo_p(y,\ep) = \{ \Phi_y^{-1}(\pi^{u,1}_y(z)) : z\in \mt\cap W^u_\ep(y) \;, \; 
\|\hf_y^p(\Phi_y^{-1}(z))\|\leq \ep\}
\subset \hLao_y(\ep)\;$$
does not coincide with $\chBo_p(y,\ep)$, however it follows from (5.1) that
$$\frac{1}{2C}\diam(\chBo_p(y,\ep)) \leq \hBo_p(y,\ep)  \leq 2C\, \chBo_p(y,\ep)\;.$$
So, it is enough to compare diameters of sets of the form $\hBo_p(y,\ep) $.

Next, notice that
\be
\|\hf^p_y(v) \| \leq C^2 \ep \quad \forall v\in \hBo_p(y,\ep)\;.
\ee
Given $x_0\in \mt$ and $x\in W^{u,1}_\ep(x_0)$, set $E^u_1(x) = T_x(W^{u,1}_\ep(x_0))$, and notice that
$d\hf^{-1}_x(0) \cdot E^{u}_1(x) = E^{u}_1(f^{-1}(x))$. 

Using the arguments in section 3 (and the proof of Lemma 3.3 above) one derives the following:

\bs 

\noindent
{\sc Lemma 6.1}
{\it Choosing $\ep_2 \in (0,\ep_1/2]$ sufficiently small, for any $x_0 \in \mt$ and any 
$x\in W^{u,1}_{\ep_2}(x_0)$ we have the following:}

\ms

(a) {\it For every $u\in E^u_1(x; \ep_2)$ there exists
$\di F_x(u) = \lim_{p\to\infty} d\hf^p_{f^{-p}(x)}(0)\cdot \hf_x^{-p}(u) \in E^u_1(x; 2\ep_2)\;.$
Moreover,  there exists a constant $C_1 > 0$ such that
$\|F_x(u) - d\hf^p_{f^{-p}(x)}(0)\cdot \hf_x^{-p}(u)\| \leq C_1\, \gamma^p\, \|u\|^2$ 
for any $u\in E^u_1(x,\ep_2)$ and any integer $p \geq 0$.}

\ms 

(b) {\it The maps $F_x : E^u_1( x;\ep_2) \longrightarrow F_x (E^u_1( x;\ep_2)) \subset E^u_1(x;2\ep_2)$ 
are $C^1$ diffeomorphisms with uniformly bounded derivatives. }

\ms

(c) {\it For  any integer $q\geq 1$  we have $d\hf_{x}^{-q}(0) \circ F_{x} (v) = F_{ f^{-q}(x)} 
\circ \hf_{x}^{-q} (v)$ for any  $v \in E^u_1 (x;\ep_2)$.}

\ms

(d) {\it For any  $\xi, u \in E^u_1(x;\ep_2/2)$ there exist the limits
$\di \Lxx  = \lim_{p\to\infty} d\hf^p_{f^{-p}(x)}(\hf_x^{-p}(\xi))\circ d\hf_x^{-p}(0)$ and
$\di F_{x,\xi}(u) = \lim_{p\to\infty} d\hf^p_{f^{-p}(x)}(\hf_x^{-p}(\xi))\circ \hf_x^{-p}(u)\;.$
Moreover, for the linear map $\Lxx : E^u_1(x) \longrightarrow E^u_1(x)$ we have  $\|\Lxx\| \leq 2$, and 
$F_{x,\xi}(u) = \Lxx \circ F_x (u)$. }

\ms

(e)   {\it For any $t \geq 0$ and any $u\in E^u_1(x;\ep_2)$ we have
$F_x(u) = \lim_{t\to \infty} d\phi_t(\phi_{-t}(x))\cdot \hphi_{x,-t}(u)$
and $d\phi_{-t}(x) \cdot F_x (u) = F_{\phi_{-t}(x)} (  \hphi_{x,-t} (u))$ ,
where $ \hphi_{x,t} = (\exp^u_{\phi_t(x)})^{-1} \circ \phi_t\circ \exp^u_x$. } \endofproof

\ms

We omit the proof, since it is almost an one-to-one repetition of the proofs of Lemma 3.3 and Theorem 3.1.

Set $\tLao_x = F_x(\hLao_x(\ep_2)) \subset E^u_1(x; 2\ep_2)$ for any $x\in \mt$.
Then, using (6.1) we get
\be
d\hf^{-1}_x(0) (\tLao_x(\ep)) \subset \tLao_{f^{-1}(x)}(\ep)\;,
\ee
and more generally $d\phi_t(x)(\tLao_x(\ep)) \subset \tLao_{\phi_t(x)}(\ep)$ for any $t \leq 0$ 
and $\ep \in (0,\ep_2]$.

It follows from (LUPC) that the distribution $E^s(x)\oplus E^0(x)\oplus E^u_1(x)$ is
integrable (see e.g. \cite{kn:Pes}), so assuming $\ep_1> 0$ is sufficiently small, 
there exist a family of invariant $C^2$ manifolds
$W^{sc}_{\ep_1}(x)$, $x\in \mt$, tangent to this distribution. 

Next, recall the local stable holonomy maps 
$\hhs_{x,y} :  W^u_{\ep_1}(x)\cap \mt  \longrightarrow W^u_{\ep_0}(y)\cap \mt$ 
($y \in \mt\cap W^s_{\ep_0}(x)$) from section 2. 
Unlike the case considered in sections 3 and 4, here there is no natural 
way to define a continuous map\footnote{Notice that the map 
$(\Phi_y)^{-1}\circ \hhs_{x,y} \circ \Phi_x$ does not necessarily send  
$\hLao_x (\ep_2)$ into $\hLao_y (\ep_1)$, since in general 
$\hhs_{x,y}$ does not map $W^{u,2}$ leaves into $W^{u,2}$ leaves.}
from $ \tLao_x (\ep_2)$ into $ \tLao_y (\ep_1)$. However, we have the following simple 
lemma which is enough to use the arguments from section 4 in the present situation.

\bs

\noindent
{\sc Lemma 6.2.}  {\it  Assuming $\ep_2 \in (0,\ep_1/2]$ is sufficiently small, for every 
$\delta\in (0,\ep_2]$, 
there exists $\delta' \in (0,\ep]$ such that for any  $y \in W^s_{\delta'}(x)\cap \mt$ and any 
$u\in \tLao_x (\ep_2)$ there exists $v\in \tLao_y(\ep_1)$ 
with $\dist (u,v) < \delta$, where $\dist$ is the distance on $TM$ induced by the Riemann metric. }

\bs

\noindent
{\it Proof of Lemma} 6.2. It is enough to deal with elements of   
$\Lao_x(\ep_2) = \pi^{u,1}_x (W^{u}_\ep(x) \cap \mt)$ and $\Lao_y(\ep_0)$.

Given $x\in \mt$, 
let $\pi_x : B(x,\ep_1) \cap \mt \longrightarrow W^{sc}_{\ep_0}(x)$  be the projection along 
leaves of $W^{u,2}$. It is well-known 
(see e.g. \cite{kn:Pes} or \cite{kn:HPS}) that $\pi_x$ is uniformly (H\"older) continuous, so given 
$\delta > 0$, there exists $\delta'' > 0$ 
such that if $z'\in W^{sc}_{\delta''}(z)\cap \mt$ for some $z\in \mt$, then 
$d(\pi_x(z), \pi_x(z')) < \delta$. Now take $\delta' > 0$ so small that if 
$y \in W^s_{\delta'}(x)\cap \mt$, then 
$d(z, \hhs_{x,y}(z)) < \delta''$ for any $z \in W^u_{\ep_1}(x)\cap \mt$. 

With this choice of $\delta'$, let $y \in W^s_{\delta'}(x)\cap \mt$. Given any $x'\in \Lao_x (\ep_2)$, 
we will show that there exists
$y' \in \Lao_y(\ep_1)$ with $d(x',y') < \delta$. Indeed,  there exists $\xi\in W^u_{\ep_1}(x)\cap \mt$
with $x' = \pi^{u,1}_x(\xi)$. Setting $\eta = \hhs_{x,y}(\xi)$ and $y' = \pi^{u,1}_y(\eta)$, we get 
$\eta \in  W^u_{\ep_1}(y)\cap \mt$,
so $y' \in \Lao_y(\ep_1)$. Moreover, $d(\xi,\eta) < \delta''$, so 
$d(x',y') = d(\pi(\xi), \pi(\eta)) < \delta$.
\endofproof

\bs

For $x\in \mt$ and $y \in \mt \cap W^u_{\ep_2}(x)$ (with a global sufficiently small constant 
$\ep_2 > 0$, as always), let
$\pi_x^y$ be the {\it projection} along $W^{u,2}$ leaves from $W^{u,1}_{\ep_2}(x)$ to 
$W^{u,1}_{\ep_1}(y)$. Initially,
$\pi_x^y$ is only defined on $\{ \pi^{u,1}_x(z) : z \in \mt \cap W^u_{\ep_0}(x)\}$. 
Moreover, the maps $\pi_x^y$ are (uniformly)
$C^1$ (see Theorem 6.1 in \cite{kn:HPS} or \cite{kn:Pes}), so taking $\ep_2 > 0$ sufficiently 
small and using Whitney's extension theorem,
we can assume that $\pi_x^y$ has a $C^1$ extension 
$\pi_x^y :  W^{u,1}_{\ep_2}(x) \longrightarrow W^{u,1}_{\ep_1}(y)$. 
Then, assuming that $\ep_3 > 0$ is sufficiently small and $y \in \mt \cap W^u_{\ep_3}(x)$, define the map
$\hpi_x^y : E^u_1(x;\ep_2) \longrightarrow E^u_1(y;\ep_1)$
by $\hpi_x^y = (\Phi_y)^{-1} \circ \pi_x^y\circ \Phi_x$.
We can assume that the constant $C > 0$ is taken so large that 
$\|\hpi_x^y(u)- \hpi_x^y(v)\| \leq C\, \|u-v\|$ for 
all $u,v\in E^u_1(x;\ep_2)$ and all $x,y$ as above.

In order to prove an analogue of Lemma 4.4 in the present situation we need the following 
which is the analogue of formula (3.14).

\bs

\noindent
{\sc Lemma 6.3.}  {\it Let $x\in \mt$,  $y \in \mt \cap W^u_{\ep_3}(x)$ and let $\eta = \hpi_x^y(0)$. 
Then for any $u \in E^u_1(x; \ep_2)$ we have}
$F_x(u) = d\hpi_x^y(\eta) \circ L_{y,\eta}\cdot \left[ F_y (\tpi_x^y(u)) - F_y(\eta)\right]\;.$
\endofproof

\ms

The proof is essentially a repetition of the proof of formula (3.14) with small modifications, 
so we omit it. 

Next, given $x\in \mt$, let $m_x \geq 1$ be the minimal integer such that there exists 
$\ep = \ep(x) \leq \ep_3$
with $\dim (\span (\tLao_x(\delta))) = m_x$ for all $0 < \delta \leq \ep$. Then the linear subspace
$E^{u,1}_\mt(x) = \span (\tLao_x(\delta))$ of $E^u_1(x)$ is the same for all $\delta \in (0,\ep]$. 
As in section 4, $m_x$ is $\phi_t$-invariant, and using the argument from the proof of Lemma 4.4 with
minor modifications, we get the following.

\bs

\noindent
{\sc Lemma 6.4.} {\it There exists an integer $m$ such that $m_x = m$ for any $x \in \mt$. Moreover,
we can choose $\ep_3 > 0$ so that for any $x\in \mt$ we have $E^{u,1}_\mt (x) = \span (\tLao_x(\ep_3))$ 
and there exists  $\ep = \ep(x)\in (0,\ep_3]$ such that $E^{u,1}_\mt(y)$ depends continuously on 
$y \in W^s_\ep(x) \cap \mt$.}
\endofproof

\bs

For $z\in \mt$, $\ep \in (0,\ep_2]$ and an integer $p \geq 0$ set
$\tB^{u,1}_p(z,\ep) = F_z(\hB^{u,1}_p(z,\ep)) \subset \tLao_z(\ep)\;.$

As in the proof of Lemma 4.1, to prove part (a) of Lemma 5.3 we have to establish the following.

\bs

\noindent
{\sc Lemma 6.5.} {\it For any  $0 < \delta \leq  \ep \leq \ep_3/2$ there exist  a  constant 
$D =  D (x,\delta , \ep) > 0$ 
and an open neighbourhood $V_0$ of $x$ in $W^s_{\ep_0}(x) \cap \mt$ such that
$\ell \left( \tB^u_p( y_p,\ep)  \right) \leq D  \, \ell \left(\tB^u_p ( y_p , \delta) \right) \;$
for any $y\in V_0$ and any integer $p \geq 0$.}

\ms

\noindent
{\it Proof of Lemma 6.5.} Choose $\ep = \ep(x) \in (0,\ep_3]$ so that $E^{u,1}_\mt(y)$
depends continuously on $y \in W^s_\ep(x)\cap \mt$. For any $y \in W^s_{\ep}(x) \cap \mt$
choose and fix an orthonormal basis $e_1(y) ,e_2(y) , \ldots,e_m(y)$ in $E^{u,1}_\mt (y)$
which depends continuously on $y$. 

Let $0 < \delta \leq  \ep \leq \ep_3/(4C^2)$. 
By the definition of $E^{u,1}_\mt(x)$ and Lemma 6.4, there exist 
$u_1, u_2, \ldots, u_m \in \tLao_x(\delta/(2C^2))$ which are 
linearly independent. Set
$\Delta = 
\Vol_m[u_1,u_2, \ldots,u_m] > 0\;,$
where $[u_1,u_2, \ldots,u_m] $ denotes the parallelepiped in $E^{u,1}_\mt(x)$ determined by the vectors
$u_1, \ldots,u_m$.
Using Lemma 6.2, choose an  open neighbourhood $V_0$ of $x$ in $W^s_{\ep}(x) \cap \mt$ such that for any
$y \in W^s_{\ep}(x) \cap \mt$ there exist $u_1(y), \ldots, u_m(y)\in \tLao_y(\delta/C^2)$ with
$\Vol_m[u_1(y),u_2(y), \ldots,u_m(y)] \geq \frac{\Delta}{2}$ 
and $\frac{\|u_j\|}{2} \leq \| u_j(y)\| \leq 2\|u_j\|$ for all $y\in V_0$, $1\leq j\leq m$. Fix such $u_j(y)$
for any $y \in V_0$ and let
$L_y = L(x,y,\delta) :  E^{u,1}_\mt(y) \longrightarrow E^{u,1}_\mt(y)$ be the linear operator 
such that  $L_y u_j(y) = e_j(y)$ 
for all $j = 1, \ldots,m$. It then follows that there exists a constant
$b = b(x,\delta) > 0$ (determined by $\Delta$ and $\| u_1\|, \ldots, \|u_m\|$) such that 
$\|L_y\| \leq b$ for all $y\in V_0$. 
 
Fix for a moment $y \in V_0$. Consider an arbitrary integer $p \geq 1$ and set $z = f^{-p}(y) \in \mt$. Given
$v \in  \tBo_p( z , \ep)$, we have $v = F_z(w)$ for some $w \in \hBo(z,\ep)$, and it follows from 
Lemma 6.1(c) and (6.2) that 
$\|d\hf_z^p(0)\cdot v\| = \|d\hf^p_z(0)\cdot F_z(w)\| 
= \| F_y(\hf_z^p(w))\|\leq 2\|\hf^p_z(w)\| \leq 2C^2\ep \leq \ep_3/2\;.$
Now $ \tBo_p( z , \ep) \subset \tLao_z(\ep)$
implies $v \in \tLao_z(\ep)$, so $u = d\hf_z^p(0)\cdot v \in \tLao_y(2C^2\ep) \subset E^{u,1}_\mt(y)$. 
Consequently, $u = \sum_{s=1}^m c_s\, u_s(y)$ for some real numbers $c_s$, so
$L_y u = \sum_{s=1}^m c_s L_y (u_s(y))  = \sum_{s=1}^m c_s e_s(y) \;.$
Thus,  $\sqrt{\sum_{s=1}^m c_s^2} = \|L_y u\| \leq \|L_y \|\, \|u\| \leq 2C^2\ep\, b \;,$ and so  
$|c_s| \leq 2C^2\ep\, b$ for all $ s = 1, \ldots, m$. 

By (6.3), $v_j = d\hf_y^{-p}(0) \cdot u_j(y) \in \tLao_{z}(\delta)$. Moreover, we have 
$v_j\in \tB^u_p(z,\delta)$ for all $j = 1,  \ldots,m$. Indeed, $u_j(y) = F_y(u'_j)$ 
for some $u'_j \in \hLao_y(\delta/C^2)$, so
$\Phi_y(u'_j) = \pi^{u,1}_y(\eta_j)$ for some $\eta_j\in \mt\cap W^u_{\delta/C}(y)$. 
Then for $v'_j = \hf^{-p}_y(u'_j)$ and $\zeta_j = f^{-p}(\eta_j) \in \mt \cap W^u_{\delta/C}(z)$ we have 
$v'_j = \Phi_z^{-1}\circ f^{-p}\Phi_y(u'_j) 
\in \hLao_z(\delta)\;,$
and therefore
$v_j =  d\hf_y^{-p}(0) \cdot F_y(u'_j) = F_z(\hf^{-p}_y(u'_j))
= F_z(v'_j) \in \tLao_z(\delta)\;.$
It follows  from $\eta_j\in W^u_{\delta/C}(y)$ that $\|\Phi^{-1}_y(\eta_j)\| \leq \delta$, so
$\|\hf_z^p(\Phi_z^{-1}(\zeta_j))\| = \|\Phi_y^{-1}(\eta_j)\| \leq \delta$ and therefore 
$v'_j \in \hBo_p(z,\delta)$. This gives $v_j = F_z(v'_j) \in \tBo_p(z,\delta)$, so
$\|v\| = \|d \hf^{-p}_y (0) \cdot u\| = 
\left\| \sum_{s=1}^m c_s\,  d\hf^{-p}_y(0) \cdot u_s(y)\right\|
 \leq  m \, 2C^2\ep \, b \, \max_{1\leq s\leq m} \|v_s\| 
 \leq  m \, 2C^2 \ep \, b \, \ell(\tB^u_p(z,\delta))$.
Hence $\ell \left( \tB^u_p( z ,\ep)  \right)  \leq D \, \ell \left( \tB^u_p ( z , \delta) \right)$,
where $D = D(x,\delta, \ep) = m\, 2C^2 \ep \, b$.
\endofproof

\bs

As in section 4, the proof of part (b) of Lemma 5.3 is essentially a repetition of the above argument,
so we omit the details.
\endofproof

\bs

\footnotesize

\noindent
{\it Acknowledgements.} Thanks are due to Keith Burns, Boris Hasselblatt, Charles Pugh and 
Amie Wilkinson for various kind of information they have provided to me.

\bs

\end{document}